\documentclass{amsart}
\usepackage[all,cmtip,ps]{xy}
\usepackage{amsmath}
\usepackage{amsthm}
\usepackage{amssymb}
\usepackage{amscd}
\usepackage[usenames]{color}

\newcommand{\A}{\mathcal{A}}
\newcommand{\B}{\mathcal{B}}
\newcommand{\Z}{\mathbb{Z}}
\newcommand*{\D}{\mathcal{D}}

\DeclareMathOperator{\Le}{\mathbf{L}}

\DeclareMathOperator{\Coker}{Coker}

\DeclareMathOperator{\Hom}{Hom}
\DeclareMathOperator{\End}{End}
\DeclareMathOperator{\Ext}{Ext}

\DeclareMathOperator{\Ker}{Ker}

\DeclareMathOperator{\soc}{soc}
\DeclareMathOperator{\rad}{rad}

\DeclareMathOperator{\imm}{Im}
\DeclareMathOperator{\Imm}{Im}

\DeclareMathOperator{\id}{id}

\newcommand*{\K}{\mathcal{K}}
\newcommand*{\Hm}{H}

\DeclareMathOperator{\injdim}{injdim}

\newcommand*{\T}{\mathcal{T}}
\newcommand*{\F}{\mathcal{F}}

\DeclareMathOperator{\R}{\mathbf{R}}
\DeclareMathOperator{\G}{\mathbf{\Gamma}}
\newcommand*{\C}{\mathcal{C}}

\newcommand*{\lMod}{\textrm{\textup{-Mod}}}
\newcommand*{\rMod}{\textrm{\textup{Mod-}}}

\newcommand*{\lmod}{\textrm{\textup{-mod}}}
\newcommand*{\rmod}{\textrm{\textup{mod-}}}
\newtheorem{thm}{Theorem}[section]

\newtheorem{lemma}[thm]{Lemma}
\newtheorem{prop}[thm]{Proposition}
\newtheorem{cor}[thm]{Corollary}
\theoremstyle{definition}
\newtheorem{DEF}[thm]{Definition}
\newtheorem{de}[thm]{Definition}

\newtheorem{rem}[thm]{Remark}

\newtheorem{EX}[thm]{Example}

\def\dualita#1#2{\mathrel{
                 \mathop{\vcenter{
                 \offinterlineskip
                 \hbox to 1.2truecm{\rightarrowfill}
                 \hbox to 1.2truecm{\leftarrowfill}}}%
                 \limits_{#2}^{#1}}}
\def\dual{\mathrel{
                 \mathop{\vcenter{
                 \offinterlineskip
                 \hbox to .6truecm{\rightarrowfill}
                 \hbox to .6truecm{\leftarrowfill}}}%
                 }}

\begin{document}
\title{Reflexivity in Derived Categories}
\author[F. Mantese]{Francesca Mantese}
\author[A. Tonolo]{Alberto Tonolo}
\address[F. Mantese]{Dipartimento di Informatica, Universit\`a degli Studi di Verona, strada Le Grazie  15, I-37134 Verona - Italy}
\email{francesca.mantese@univr.it}
\address[A. Tonolo]{ Dip. Matematica Pura ed Applicata, Universit\`a degli studi di Padova, via Trieste 63, I-35121 Padova Italy}
\email{tonolo@math.unipd.it}
\thanks{Research supported by grant CPDA071244/07 of Padova University}
\dedicatory{to the memory of our friend and colleague Silvia Lucido}
\begin{abstract}
An adjoint pair of contravariant functors between abelian categories can be extended to the adjoint pair of their derived functors in the associated derived categories. We describe the reflexive complexes and interpret the achieved results in terms of objects of the initial abelian categories. In particular we prove that, for functors of any finite cohomological dimension, the objects of the initial abelian categories which are reflexive as stalk complexes form the largest class where a {\color{black}Cotilting Theorem in the sense of Colby and Fuller \cite[Ch.~5]{CbF1} works}.
\end{abstract}

\maketitle

\section*{Introduction}

Adjoint pairs of functors in derived categories are deeply studied by several authors (see for instance \cite{ CPS, Hp, H, R}).
A large class of these pairs is obtained extending (\cite[Lemma~13.6]{K}) adjoint functors between abelian categories to their derived functors.
In this paper we focus on the adjunctions and on the corresponding dualities obtained by extending contravariant adjoint functors. Wide classes of examples arise in module and sheaf categories. 

Both the adjunctions, the one in the abelian categories and the one in the associated derived categories, determine a different notion of reflexivity: to distinguish them, we will call
\emph{$\D$-reflexive} the complexes which are reflexive with respect to the adjunction in the derived categories. An object of the starting abelian categories is called \emph{$\D$-reflexive} if it is $\D$-reflexive as stalk complex, simply \emph{reflexive} if it is reflexive with respect to the adjunction in the abelian categories.

Our main aim is to describe the $\D$-reflexive complexes and to study the $\D$-reflexive objects of the initial abelian categories. Reading on the underlying abelian categories all the effects of the duality in the corresponding derived categories,
we prove that for an adjoint pair of functors of cohomological dimension 1, the $\D$-reflexive objects are exactly those for which a {\color{black}\emph{Cotilting Theorem} in the sense of Colby and Fuller (\cite[Ch.~5]{CbF1}) holds (Theorem~\ref{cor:BB}).}
Our approach allows on one side to read this celebrated result in its traditional framework as a natural consequence of the duality between derived categories induced by contravariant $\Hom$-functors associated to a cotilting bimodule. On the other side it permits to generalize the Cotilting Theorem to arbitrary abelian categories and  adjoint pairs of contravariant functors of any finite cohomological dimension. Thus, on the one hand we get a general and unitary version of all the several cases considered in the literature of dualities induced by cotilting bimodules of injective dimension 1 (see  \cite{Cb, Cb1, CbF, CbF1, C, CF, Ma, T}). On the other hand, under cohomological conditions automatically satisfied in the traditional settings, we succeed in finding positive results which generalize to arbitrary abelian categories and adjoint functors of cohomological dimension $n$ the results obtained by Miyashita \cite{M} for a cotilting bimodule of injective dimension $n$ in the noetherian case.


In the first section we recall some preliminaries on derived functors and their properties; particular attention is dedicated to the notion of \emph{way-out} functor.

In the second section we describe the adjoint pairs we are interested in, and we compare  the related notions of reflexivity  in the abelian and in the associated derived categories. In particular we give examples of $\D$-reflexive objects in the starting abelian categories which are not reflexive and, conversely, of reflexive objects which are not $\D$-reflexive.

In the third section we investigate the relation between the $\D$-reflexivity of a complex and of its terms or its cohomologies. We show that if the cohomologies or the terms are $\D$-reflexive, then so is the complex itself. The converse in general is not true (see Examples~\ref{ex:Gdim}, \ref{ex:GGdim}). Assuming that the functors have cohomological dimension at most one, we prove that a complex is $\D$-reflexive if and only if its cohomologies are $\D$-reflexive (see Corollary~\ref{cor:atmostone}).

In the fourth and fifth sections we study in details the $\D$-reflexive objects in the abelian categories we start from. The fourth is dedicated to the favorable case of functors of cohomological dimension $\leq 1$. We prove that a Cotilting Theorem \cite{CbF1} for the classes of $\D$-reflexive objects (Theorem~\ref{cor:BB}) holds. 
Finally the fifth section is devoted to the case of functors with arbitrary finite  cohomological dimension, assuming the abelian categories have enough projectives. The latter hypothesis permits us to use the standard tool of spectral sequences. For a spectral sequence analysis in the covariant case see \cite{BB1}. This approach allows us to reveal the cohomological conditions (Condition I, II page~\pageref{conditionI}) necessary to generalize the results obtained in cohomological dimension $\leq 1$ (Theorems~\ref{thm:last}, \ref{thm:lastt}). In particular we completely recover the results obtained for module categories by Miyashita \cite{M} in the noetherian case and in \cite{AT} for arbitrary associative rings.

Several examples occur along the whole paper describing pathologies and positive results.


For the unexplained notations in module theory we refer to \cite{AF}, for those in sheaf theory to \cite{H^1}. 
We follow  \cite{H, W} for definitions and results regarding derived categories, derived functors and  spectral sequences.

\section{Preliminaries}

Given an abelian category $\mathcal A$, we denote by $\K(\mathcal A)$ (resp. $\K^+(\mathcal A)$, $\K^-(\mathcal A)$, $\K^b(\mathcal A)$) the homotopy category of unbounded (resp. bounded below, bounded above, bounded) complexes of objects of $\A$ and by $\D(\mathcal A)$ (resp. $\D^+(\mathcal A)$, $\D^-(\mathcal A)$, $\D^b(\mathcal A)$) the associated derived category. 
In the sequel with $\D^*(\mathcal A)$ or $\D^{\dag}(\mathcal A)$ we will denote any of these derived categories. Moreover, with $\D^*_{\leq n}(\mathcal A)$, $n\in\mathbb Z$, we mean all the complexes in $\D^*(\mathcal A)$  whose cohomologies are zero in any degree greater than $n$. Similarly, we define $\D^*_{\geq n}(\mathcal A)$.

All the considered functors between derived categories are assumed to be \emph{$\delta$-functors}, i.e. they commute with the shift functor and send triangles to triangles.
 Given an object $M\in\mathcal A$, we  continue to denote by $M$ also the \emph{stalk complex} in $\D(\mathcal A)$ associated to $M$, i.e. the complex with $M$ concentrated in degree zero.

Let $X: \dots \to X_{-1}\stackrel{g_{-1}}{ \to} X_0 \stackrel{g_0}{\to} X_1 \to \dots$ be a complex in $\D(\mathcal A)$. For any integer $n\in \Z$ we define the following truncations:
\[ \tau_{\scriptscriptstyle >n}X: \dots\to 0 \to X_{n+1}\to X_{n+2}\to \dots \quad   \tau_{\scriptscriptstyle \leq n} X: \dots \to X_{n-1}\to X_n\to 0\to\dots\]
\[ \sigma_{\scriptscriptstyle >n} X :    \dots\to 0\to  X_n/\ker{g_n} \to X_{n+1}\to \dots \quad   \sigma_{\scriptscriptstyle \leq n}X:  \dots\to X_{n-1}\to \ker{g_n}\to 0\to\dots\]
 In particular, for any $n\in\Z$ there are the following triangles:
\[ \tau_{\scriptscriptstyle > n} X \to X \to \tau_{\scriptscriptstyle \leq n}X \to  
 \tau_{\scriptscriptstyle >n}X[1]\qquad \sigma_{\scriptscriptstyle \leq n} X\stackrel{}{\to} X\stackrel{\pi_X}{\to} \sigma_{\scriptscriptstyle >n} X \to \sigma_{\scriptscriptstyle \leq n} X[1].\]
 
 In this section we study the behavior of the composition of contravariant \emph{way-out functors} and the relations among the way-out conditions, the finite cohomological dimension and the closure properties of the acyclic objects associated to a contravariant functor.
 Let us first recall the definition of \emph{way-out} functors, as in \cite[Chp. I \S 7]{H} and \cite[Chp. I \S 11]{L}.

\begin{de}
Let $\A$ and $\B$ be abelian categories and let $F: \D^*(\A) \to \D(\B)$ be a covariant (resp. contravariant) functor.
\begin{enumerate}
\item The functor $F$ is \emph{way-out left} if there exists $n\in\mathbb Z$ such that
\[F(\D_{\leq 0}^*(\A))\subseteq \D_{\leq n}(\B)\quad(\text{resp. }
F(\D_{\geq 0}^*(\A))\subseteq \D_{\leq n}(\B));
\]
in such a case we define the \emph{upper dimension} of $F$ setting
\[\dim^+ F=\inf\{n:F(\D_{\leq 0}^*(\A))\subseteq \D_{\leq n}(\B)\}
\quad
(\text{resp. } =\inf\{n:F(\D_{\geq 0}^*(\A))\subseteq \D_{\leq n}(\B)\}).\]
\item The functor $F$ is \emph{way-out right} if there exists $n\in\mathbb Z$ such that
\[F(\D_{\geq 0}^*(\A))\subseteq \D_{\geq n}(\B)\quad(\text{resp. }
F(\D_{\leq 0}^*(\A))\subseteq \D_{\geq n}(\B))
;
\]
in such a case we define the \emph{lower dimension} of $F$ setting
\[\dim^- F=\sup\{n:F(\D_{\geq 0}^*(\A))\subseteq \D_{\geq n}(\B)\}\quad
(\text{resp. } =\sup\{n:F(\D_{\leq 0}^*(\A))\subseteq \D_{\geq n}(\B)\}).\]
\end{enumerate}
\end{de}

\begin{rem}
Let $F: \D^*(\A) \to \D(\B)$ be a covariant (resp. contravariant) functor. If $F$ is way-out left and $\dim^+F=m$, then for any $k\in\Z$
\[F(\D_{\leq k}^*(\A))\subseteq \D_{\leq k+m}(\B)\quad(\text{resp. }
F(\D_{\geq k}^*(\A))\subseteq \D_{\leq m-k}(\B))\]
Analogously, if $F$ is way-out right and $\dim^-F=m$, then for any $k\in\Z$
\[F(\D_{\geq k}^*(\A))\subseteq \D_{\geq k+m}(\B)\quad(\text{resp. }
F(\D_{\leq k}^*(\A))\subseteq \D_{\geq m-k}(\B)).\]
Clearly, if $F$ is both way-out left and right, then it is \emph{bounded}, i.e. it sends  bounded complexes in $\D^*(\A)$ to bounded complexes in $\D(\B)$.
\end{rem}

The following easy proposition will be useful in the sequel. 

\begin{prop}\label{prop:way}
Let $G_1: \D^*(\A) \to \D^{\dag}(\B)$ and $G_2: \D^{\dag}(\B) \to \D(\C)$ be two contravariant functors and $F=G_2 G_1$ their composition.
\begin{enumerate}
\item If $G_1$ is way-out left with $\dim^+G_1=m_1$ and $G_2$ is way-out right with $\dim^-G_2=m_2$, then $F$ is way-out right with $\dim^-F=m_2-m_1$;
\item if $G_1$ is way-out right with $\dim^-G_1=m_1$ and $G_2$ is way-out left with $\dim^+G_2=m_2$, then $F$ is way-out left with $\dim^+F=m_2-m_1$.
\end{enumerate}
\end{prop}

From now on we denote by $\Phi:\mathcal A\to\mathcal B$ and $\Psi:\mathcal B\to \mathcal A$ two additive non zero contravariant functors between the abelian categories $\A$ and $\B$. Following \cite[Theorem~5.1]{H}, to guarantee the existence of the derived functors $\R^*\Phi:\D^*(\mathcal A)\to \D(\mathcal B)$ and $\R^{\dag}\Psi:\D^{\dag}(\mathcal B)\to \D(\mathcal A)$, we assume the existence of triangulated subcategories $\mathcal P$ of $K^*(\mathcal A)$ and $\mathcal Q$ of  $K^{\dag}(\mathcal B)$ such that:
\begin{itemize}
\item every object of $K^*(\mathcal A)$ and every object of $K^{\dag}(\mathcal B)$ admits a quasi-isomorphism from objects of $\mathcal P$ and $\mathcal Q$, respectively;
\item if $P$ and $Q$ are exact complexes in $\mathcal P$ and $\mathcal Q$, then also $\Phi(P)$ and $\Psi(Q)$ are exact.
\end{itemize}

Given complexes $X\in\D^*(\mathcal A)$ and 
$Y\in\D^{\dag}(\mathcal B)$, we have
$\R^*\Phi X= \Phi P$ and $\R^{\dag}\Psi Y= \Psi Q$, where $P$ is a complex in $\mathcal P$ quasi-isomorphic to $X$, and $Q$ is a complex in $\mathcal Q$ quasi-isomorphic to $Y$. 

If $\Phi(K^*(\A))\subseteq K^{\dag}(\B)$ and
$\Phi(\mathcal P)\subseteq\mathcal Q$, then there exists also $\R^*(\Psi\Phi)$ and it is isomorphic to $\R^{\dag}\Psi\R^*\Phi$ \cite[Proposition~5.4]{H}.

 \begin{de}
 \begin{enumerate}
 \item An object $A$ in $\mathcal A$ is called \emph{$\Phi$-acyclic} if $H^i(\R^*\Phi A)=0$ for any $i\not=0$.
 \item The category $\mathcal A$ has \emph{enough $\Phi$-acyclic objects}  if any object in $\mathcal A$ is image of a $\Phi$-acyclic object.
 \item The functor $\Phi$ has \emph{cohomological dimension $\leq n$} if, for each $A$ in $\mathcal A$, we have $H^i(\R^*\Phi A)=0$  for $|i|>n$
 \end{enumerate}
 \end{de}
 
 \begin{rem}
If $\mathcal A$ has enough $\Phi$-acyclic objects, then the right derived functor $\R^- \Phi:\D^-(\mathcal A)\to \D({\mathcal B})$ is defined and it may be computed using $\Phi$-acyclic resolutions: given a complex $X\in\D_{\leq n}(\mathcal A)$, we have that  $\R^-\Phi X= \Phi L$, where $L$ is a complex in $\D_{\leq n}(\mathcal A)$ with $\Phi$-acyclic terms quasi-isomorphic to $X$. 
In particular, if the category $\A$ has enough projectives, $\R^- \Phi:\D^-(\mathcal A)\to \D({\mathcal B})$ is defined, and for each object $A$ in $\mathcal A$, $H^n(\R^-\Phi A)$ coincides with the usual right $n^{th}$-derived functor of $\Phi$ evaluated in $A$.

If $\Phi$ has finite cohomological dimension $n$ and $\A$ has enough $\Phi$-acyclics, then any complex $X\in\D(\A)$ is quasi-isomorphic to a complex $L$ with $\Phi$-acyclic terms; thus the total derived functor $\R\Phi$ exists and $\R\Phi X=\Phi L$ (see \cite[Corollary~5.3, $\gamma$.]{H}).
\end{rem}

Notice that if $\R^*\Phi:\D^*(\mathcal A)\to \D({\mathcal B})$ is way-out in both directions, then it has finite cohomological dimension. Under the hypothesis that $\A$ has enough $\Phi$-acyclic objects, also the converse holds:
\begin{prop}\label{prop:rightleft}
Let $\A$ and $\B$ be abelian categories, and $\Phi:\A \to \B$ a contravariant functor. Assume $\A$ has enough $\Phi$-acyclic objects; then
\begin{enumerate}
\item $\R^-\Phi$ is way-out right of lower dimension $ 0$;
\item if $\Phi$ has finite cohomological dimension $n$,  then $\R\Phi$ is way-out left of upper dimension $n$.
\end{enumerate}
\end{prop}
\begin{proof}
1. It is clear, since $\R^-\Phi$ may be computed on $\Phi$-acyclic resolutions.

2. Let $X:=0\to X_0\to X_1\to ...$ be an object in $\D_{\geq 0}(\A)$. Since $\Phi$ has finite cohomological dimension, there exists a complex of $\Phi$-acyclic objects $L:=...\to L_{-1}\to L_0\to L_1\to ...$ quasi isomorphic to $X$. Denote by $C$ the cokernel of $L_{-1}\to L_0$; then $C$ is quasi isomorphic to
$\tau_{\scriptscriptstyle\leq 0}L$.
Since $H^i(\R\Phi C)=0$ for $i>n$, then for each $i>n$ we have
\[0=H^i(\R\Phi(\tau_{\scriptscriptstyle\leq 0}L))=H^i(\Phi(\tau_{\scriptscriptstyle\leq 0}L))=H^i(\Phi L)=H^i(\R\Phi X ).\] 
\end{proof}

\begin{prop}\label{lemma:aciclici}
 Assume $\R^*\Phi$ is way-out right of lower dimension $\geq 0$ and way-out left of upper dimension $\leq n$. If $0\to X_{0}\to X_{1} \to \dots \to X_n$ is an exact complex where the $X_i$, $i>0$, are $\Phi$-acyclic objects of $\mathcal A$, then also $X_{0}$ is $\Phi$-acyclic. In particular, if $n=1$ the class of $\Phi$-acyclic objects is closed under submodules.
\end{prop}
\begin{proof}
Let $X:=0\to X_{0}\to X_{1} \to \dots \to X_n\to 0$; we have to prove that $H^i(\R^*\Phi(X_{0}))=0$ for each $i\not = 0$. Since the stalk complex $X_{0}$ belongs to $\D^*_{\leq 0}(\mathcal A)$, and $\R^*\Phi$ has lower dimension $\geq 0$, $\R^*\Phi(X_{0})$ belongs to $\D_{\geq 0}(\mathcal B)$. Therefore it is sufficient to prove that $H^i(\R^*\Phi(X_{0}))=0$ for each $i> 0$.
Consider the triangle
\[\tau_{\scriptscriptstyle>0}X\to X\to \tau_{\scriptscriptstyle\leq 0}X\to
\tau_{\scriptscriptstyle>0}X[1]\]
and observe that $\tau_{\scriptscriptstyle\leq 0}X$ is the stalk complex $X_0$; then for each $i>0$ we have the exact sequence 
\[ H^{i-1}(\R^*\Phi(\tau_{\scriptscriptstyle>0}X))\to H^{i}(\R^*\Phi(X_0))\to
H^{i}(\R^*\Phi(X)).\]
Now, since the terms in $\tau_{\scriptscriptstyle>0}X$ are $\Phi$-acyclics, $\R^*\Phi(\tau_{\scriptscriptstyle>0}X)$ has non-zero terms only in negative degrees, 
and therefore it has the $(i-1)^{th}$ cohomology equal to zero. Since $X$ belongs to
$\D^*_{\geq n}(\mathcal A)$, and $\dim^+(\R^*\Phi)\leq n$,  we get that $\R^*\Phi(X)$ belongs to $\D^*_{\leq 0}(\mathcal B)$. So also the $i^{th}$ cohomology of $\R^*\Phi(X)$ vanishes, and we conclude.
\end{proof}

 \begin{de}\label{de:aciclici}
 Assume that $\Phi(K^*(\A))\subseteq K^{\dag}(\B)$. We say that an object $L\in\A$ is \emph{$\Psi$-$\Phi$-acyclic} if $L$ is $\Phi$-acyclic and $\Phi(L)$ is $\Psi$-acyclic. We say that the abelian category $\A$ has \emph{enough $\Psi$-$\Phi$-acyclic objects} if any $A\in \A$ is image of a $\Psi$-$\Phi$-acyclic object.
\end{de}

%
 
\begin{prop}\label{rem:dim}
Assume that $\Phi(K^*(\A))\subseteq K^{\dag}(\B)$.
\begin{enumerate}
\item If $\A$ has enough $\Psi$-$\Phi$-acyclic objects, then $\R^{\dag}\Psi\R^*\Phi$  is way-out left of upper dimension $\leq 0$. 
\item If $\Phi$ has cohomological  dimension $n$, and $\A$ and $\B$ have enough $\Phi$-acyclics and $\Psi$-acyclics, respectively, then $\R^{\dag}\Psi\R^*\Phi$ 
 is way-out right of lower dimension $\geq -n$.
\end{enumerate}
 \end{prop}
\begin{proof}
1. It follows since any complex in $\D^*_{\leq 0}(\A)$ is quasi-isomorphic to a complex in $\D^*_{\leq 0}(\A)$ with $\Psi$-$\Phi$-acyclic terms.  

2. It follows by Propositions~\ref{prop:way}, \ref{prop:rightleft}.
\end{proof}

\section{Adjunction and Reflexive objects}

From now on, we are interested in the situation when $(\Phi, \Psi)$ is a right adjoint pair; in particular $\Phi$ and $\Psi$ are left exact. The following result has a key role in our analysis.

\begin{lemma}[{\cite[Lemma~13.6]{K}}]\label{lemma:CPS}
Let  $(\Phi,\Psi)$ be a right adjoint pair. Assume that $\R^*\Phi(X)$ and $\R^{\dag}\Psi(Y)$ belong to $\D^{\dag}(\mathcal B)$ and $\D^*(\mathcal A)$, for any $X$ in $\D^*(\mathcal A)$ and $Y$ in $\D^{\dag}(\mathcal B)$, respectively. Then $(\R^* \Phi, \R^{\dag} \Psi)$ is a right adjoint pair.
\end{lemma}

Thus, under suitable assumptions on the existence of the derived functors, any adjunction in abelian categories can be extended to the associated derived categories. In this section we compare these two adjunctions. In particular we describe the relationship between the units of the two adjunction, and we show the independence of the related notions of reflexivity.

\begin{EX}\label{ex:fasci}
1. Let $(X, \mathcal O_X)$ be a locally noetherian scheme such that every coherent sheaf on $X$ is a quotient of a locally free sheaf. Consider the abelian category $\mathfrak{Mod} X$ of sheaves of $\mathcal O_X$-modules, and the thick subcategory $\mathfrak{Coh}X$ of coherent sheaves. Let $\mathcal G$ be a coherent sheaf of finite injective dimension; consider the functor $\mathcal{H}om(-, \mathcal G):\mathfrak{Mod} X\to \mathfrak{Mod} X$. The pair $(\mathcal{H}om(-, \mathcal G), \mathcal{H}om(-, \mathcal G))$ is a right adjunction.
By   \cite[Chp. III]{H^1} there exists the derived functor
\[\R^b\mathcal{H}om(-, \mathcal G):\D^b({\mathfrak{Coh}X})\to
\D^b({\mathfrak{Coh}X}).\]
 Therefore, by Lemma~\ref{lemma:CPS}, $(\R^b\mathcal{H}om(-, \mathcal G), \R^b\mathcal{H}om(-, \mathcal G))$ is a right adjoint pair. 

Moreover, by \cite[Cor. I.5.3]{H}, there exists also  the total derived functor
\[\R\mathcal{H}om(-, \mathcal G):\D({\mathfrak{Coh}X})\to
\D({\mathfrak{Coh}X}).\]
and so $(\R\mathcal{H}om(-, \mathcal G), \R\mathcal{H}om(-, \mathcal G))$ is a right adjoint pair.

2. Let $R$ be a ring, $_RU$ a left $R$-module and $S$
the endomorphism ring of $_RU$. The pair $(\Hom_R(-,U), \Hom_S(-, U))$ is a right adjunction. By \cite[Theorem C]{S}, the derived functors $\R\Hom_R(-,U) $ and $\R\Hom_S(-, U)$ always exist  and so $(\R\Hom_R(-,U), \R\Hom_S(-, U))$ is a right adjoint pair. If both $_RU$ and $U_S$ have finite injective dimension,  then $\R\Hom_R(-,U)$ and  $\R\Hom_S(-, U)$ are bounded, since they  are way-out in both directions. It follows that also $(\R^b\Hom_R(-,U), \R^b\Hom_S(-, U))$ is a right adjoint pair.
\end{EX}

In the sequel, we assume  that $(\Phi, \Psi)$ is an adjoint pair inducing the adjoint pair $(\R^*\Phi,\R^\dag\Psi)$.
Denoted by $\eta$ and $\xi$ the units of the right adjoint pair $(\Phi,\Psi)$, we indicate with $\hat{\eta}$ and $\hat{\xi}$ the units of the right adjoint pair $(\R^* \Phi, \R^{\dag} \Psi)$, i.e. the natural maps 
\[\hat{\eta}: id_{D^*(\mathcal A)}\to \R^{\dag} \Psi\R^* \Phi, \quad \hat{\xi}:id_{D^{\dag}(\mathcal B)}\to \R^*\Phi\R^{\dag}\Psi\]
such that $\R^* \Phi(\hat{\eta}_{X})\circ \hat{\xi}_{\R^* \Phi X}=1_{\R^* \Phi X}$ and $\R^{\dag}\Psi( \hat{\xi}_{Y})\circ\hat{\eta}_{\R^{\dag}\Psi Y}=1_{\R^{\dag}\Psi Y}$ for each $X$ in $D^*(\mathcal A)$ and each $Y$ in $D^{\dag}(\mathcal B)$.

Suppose that  $\A$ has enough $\Psi$-$\Phi$-acyclic objects. Let $X\in\D^-(\mathcal A)$ and $L$ be a complex in $\K^-(\mathcal A)$  of $\Psi$-$\Phi$-acyclics quasi-isomorphic to $X$. Then $\R^{\dag} \Psi\R^* \Phi X=\Psi\Phi L$,  
$\hat{\eta}_{X}$ is isomorphic to $\hat{\eta}_{L}$ in $\D^-(\mathcal A)$, and the latter coincides with the term to term extension of the unity $\eta$ to $\K^-(\mathcal A)$.

\begin{prop}\label{prop:legame}
Assume that  $\A$ has enough $\Psi$-$\Phi$-acyclic objects. Let $A$ be an object of $\A$ and $\iota$ be the canonical map of complexes $\sigma_{\leq 0}\R^*\Phi A\to \R^*\Phi A$. If $\Phi(A)$ admits a $\Psi$-acyclic resolution, we have
\[\eta_A=H^0(\R^{\dag}\Psi(\iota))\circ H^0(\hat\eta_A)=H^0(\R^{\dag}\Psi(\iota)\circ \hat\eta_A).\] 
\end{prop}
\begin{proof}
Consider a $\Psi$-$\Phi$-acyclic resolution $L: ...\to L_{-1}\stackrel{d_{-1}}{\to} L_0\to 0$ of $A$ with augmentation $f:L_0\to A$; we have the commutative diagram:
\[
\xymatrix{
L_0\ar[rr]^{f}\ar[d]^{\hat\eta_{L_0}=\eta_{L_0}}&&A\ar[d]^{\hat\eta_A}\\
\R^{\dag}\Psi\R^*\Phi L_0=\Psi\Phi(L_0)\ar[rr]^-{\R^{\dag}\Psi\R^*\Phi(f)}&&\R^{\dag}\Psi\R^*\Phi A\ar[rr]^-{\R^{\dag}\Psi(\iota)}&&\R^{\dag}\Psi(\sigma_{\leq 0}\R^*\Phi A)=\R^{\dag}\Psi(\Phi A)
}
\]
Applying the cohomology functor $H^0$, the solid part of the following diagram commutes:
\[
\xymatrix{
L_0\ar[rr]^{f}\ar[d]^{\hat\eta_{L_0}=\eta_{L_0}}&&A\ar@{-->}@/^13pt/[drr]^{\eta_A}\ar[d]^{H^0(\hat\eta_A)}\\
\Psi\Phi(L_0)\ar@{-->}@/_23pt/[rrrr]_{\Psi\Phi(f)}\ar[rr]^-{H^0(\R^{\dag}\Psi\R^*\Phi(f))}&&H^0(\R^{\dag}\Psi\R^*\Phi A)\ar[rr]^-{H^0(\R^{\dag}\Psi(\iota))}&&\Psi(\Phi A)
}
\]
Let us see that $H^0(\R^{\dag}\Psi(\iota))\circ H^0(\R^{\dag}\Psi\R^*\Phi(f))=\Psi\Phi(f)$; then, for the naturality of $\eta$ we will have 
\[\eta_A\circ f=\Psi\Phi(f)\circ \eta_{L_0}=(H^0(\R^{\dag}\Psi(\iota))\circ H^0(\R^{\dag}\Psi\R^*\Phi(f))\circ\hat\eta_{L_0}=
(H^0(\R^{\dag}\Psi(\iota))\circ H^0(\hat\eta_A))\circ f;\] 
since $f$ is an epimorphism, we will conclude.

Let $Q$ be a $\Psi$-acyclic resolution of $\Phi(A)$. Consider the diagram
\[
\xymatrix{
0\ar[r]&\Phi(L_0)\ar@{=}[d]\ar[r]&0&&&\R^*\Phi L_0\\
0\ar[r]\ar[u]&\Phi(L_0)\ar[r]&\Phi(L_{-1})\ar[u]\ar[r]&...&&\R^*\Phi A\ar[u]_{\R^*\Phi(f)}\\
0\ar[r]\ar[u]&\Phi(A)\ar[u]_{\Phi(f)}\ar[r]&0\ar[u]\ar[r]&...&&\Phi A\ar[u]\\
Q_{-1}\ar[r]\ar[u]&Q_0\ar[u]\ar[r]&0\ar[u]&&&Q\ar[u]^{qiso}
}
\]
Applying $\Psi$ we get the commutative diagram

\[
\xymatrix{
0\ar[d]&&\Psi\Phi(L_0)
\ar@{=}[d]
\ar[ll]&0\ar[d]\ar[l]&&\R^{\dag}\Psi\R^*\Phi (L_0)\ar[d]^{\R^{\dag}\Psi\R^*\Phi(f)}\\
0\ar[dd]&&\Psi\Phi(L_0)
\ar[dd]\ar[dl]^{\Psi\Phi(f)}\ar[ll]&\Psi\Phi(L_{-1})\ar[dd]\ar[l]_{\Psi\Phi(d_{-1})}&...\ar[l]&\R^{\dag}\Psi\R^*\Phi (A)\ar[dd]^{\R^{\dag}\Psi(\iota)}\\
&\Psi\Phi(A)\ \ \ar@{^(->}[dr]\\
\Psi(Q_{-1})&&\Psi(Q_0)\ar[ll]&0\ar[l]&&\R^{\dag}\Psi(\Phi(A))
}
\]
Therefore, having observed that $\Psi\Phi(f)\circ\Psi\Phi(d_{-1})=0$, the induced maps on the $0^{th}$-cohomologies are obtained as follows:
\[\xymatrix{
&H^0(\R^{\dag}\Psi\R^*\Phi L_0)\ar@{-->}[d]_{H^0(\R^{\dag}\Psi\R^*\Phi(f))}\ar@{=}[r]&\Psi\Phi(L_0)\ar@{=}[d]\\
\Coker(\Psi\Phi(d_{-1})=\!\!\!\!\!\!\!\!\!\!\!\!\!\!\!\!
&H^0(\R^{\dag}\Psi\R^*\Phi A)\ar@{-->}[d]_{H^0(\R^{\dag}\Psi(\iota))}&\Psi\Phi(L_0)\ar@{->>}[l]\ar[d]^{\Psi\Phi(f)}\\
&H^0(\R^{\dag}\Psi(\Phi A))\ar@{=}[r]&{}\Psi\Phi(A)
}
\]
i.e. $\Psi\Phi(f)=H^0(\R^{\dag}\Psi(\iota))\circ H^0(\R^{\dag}\Psi\R^*\Phi(f))$.
\end{proof}

Both the adjoint pairs $(\Phi, \Psi)$ and $(\R^*\Phi,\R^\dag\Psi)$ define on the corresponding categories the classes of \emph{reflexive objects}, i.e. the classes where the unity maps induce isomorphisms. To distinguish, we call simply \emph{reflexive} the objects $A$ in $\A$ or $B$ in $\B$ such that the natural maps $\eta_A$ or $\xi_B$ are isomorphisms; instead we say \emph{$\D$-reflexive} the complexes which are reflexive with respect to the adjoint pair $(\R^*\Phi,\R^\dag\Psi)$.
Observe that
any object $A$ in $\A$ is also, in a natural way, an object in $\D^*(\A)$. Both the maps $\eta_A$ and $\hat\eta_A$ can be considered; therefore $A$ can be reflexive or $\D$-reflexive. The two notions are independent:

\begin{EX}\label{ex:reflDrefl}
In this and all future examples $k$ denotes an algebraically closed field. For any finite-dimensional $k$-algebra given by a quiver with relations, if $i$ 
is a vertex, we 
denote by $P(i)$ the indecomposable projective associated to $i$, 
by $E(i)$ the indecomposable injective associated to $i$, and 
by $S(i)$ the simple top of $P(i)$ or, equivalently,
the simple socle of $E(i)$.

Let $\Lambda$ denote the $k$-algebra given by the quiver ${\cdot}1\stackrel 
a{\to} {\cdot}2\stackrel b{\to} {\cdot}3\stackrel c{\to} {\cdot}4$ with relations $ba=0=cb$. 

(1) Let ${}_{\Lambda}W=S(1)\oplus S(3)$; then $S=\End {}_{\Lambda}W$ is $k\oplus k$. Since ${}_{\Lambda}W$ and $W_S$ have finite injective dimension, we have the two right adjoint pairs
\[
(\Hom_{\Lambda}(-,W),\Hom_S(-,W))\text{ and }(\R^b\Hom_{\Lambda}(-,W),\R^b\Hom_S(-,W)).
\]
An easy computation permits to verify that  $S(1)$ is reflexive. Regarding $S(1)$ as a stalk complex, it is quasi isomorphic to its projective resolution $P:=0\to P(3)\to P(2)\to P(1)\to 0$. Since $P$ has $\Hom_S(-,W)$-$\Hom_{\Lambda}(-,W)$-acyclic terms,
\[\R^b\Hom_S(\R^b\Hom_{\Lambda}(S(1),W),W)=\Hom_S(\Hom_{\Lambda}(P,W),W)\]
is the complex $0\to S(3)\to 0\to S(1)\to 0$,
which is not quasi-isomorphic to $P$. Then $S(1)$ is not $\D$-reflexive.

(2) Let $_{\Lambda}{\Lambda}_{\Lambda}$ be the regular bimodule. Since the left and the right regular modules have finite injective dimension, we have the two right adjoint pairs
\[
(\Hom_{\Lambda}(-,\Lambda),\Hom_\Lambda(-,\Lambda))\text{ and }(\R^b\Hom_{\Lambda}(-,\Lambda),\R^b\Hom_\Lambda(-,\Lambda)).
\]
It is straightforward to verify that the simple module $S(2)\in\Lambda\lmod$ is not reflexive. Since all indecomposable projective modules are reflexive and $\Hom_{\Lambda}(-,\Lambda)$-$\Hom_\Lambda(-,\Lambda)$-acyclic, the simple module $S(2)$ is $\D$-reflexive.
\end{EX}

\begin{prop}\label{prop:refl}
If $A\in\mathcal A$ is reflexive and $\Psi$-$\Phi$-acyclic, then $A$ is $\D$-reflexive.
\end{prop}
\begin{proof}
We have to prove that $\hat\eta_A$ is a quasi-isomorphism, i.e. $H^i(\hat\eta_A)$ are isomorphisms for each $i$. Since $A$ is $\Psi$-$\Phi$-acyclic, then $\R^{\dag}\Psi\R^*\Phi A$ is the stalk complex $\Psi\Phi A$. Clearly $H^i(\hat\eta_A)=0$ for each $i\not=0$ are isomorphisms;   since $A$ is reflexive, also $H^0(\hat\eta_A)= \eta_A$ is an isomorphism.
\end{proof}

The category of $\D$-reflexive complexes is a triangulated subcategory of $\D^*(\A)$.  In particular the subcategory of stalk $\D$-reflexive complexes is thick, i.e, if two terms of a short exact sequence in $\mathcal A$ are $\D$-reflexive,  then also the third is.  This follows easily since any short  exact sequence in $\mathcal A$ gives rise to a triangle in $\D^*(\mathcal A)$. 

Note that, from the adjunction formulas, it follows that if a complex $X$ is $\D$-reflexive, then also $\R^* \Phi X$ is $\D$-reflexive.

\begin{DEF}\cite[Sect. 2]{AC}\label{def:cot}
Let $R$ be a ring.
 A left module $_RU$ is \emph{partial cotilting}  if it satisfies the following conditions: 
\begin{enumerate}
\item $\injdim {}_RU<\infty$;
\item
  $\Ext^i_{R}(U^{\alpha}, U)=0$, for each $i>0$ and any cardinal $\alpha$.
\end{enumerate}
The module $_RU$ is \emph{cotilting} if  moreover the following condition is satisfied
\begin{enumerate}
\item[(3)] there exists $n\in\mathbb N$ and an exact sequence $0\to U_n\to \dots\to U_1\to U_0\to Q\to 0$ where $Q$ is an injective cogenerator of $R\lMod$ and $U_i$ are direct summands of products of copies of $U$.
\end{enumerate}
 A bimodule $_RU_S$ is  \emph{(partial) cotilting}  if both $_RU$ and $U_S$ are (partial) cotilting.
 \end{DEF}

Partial cotilting modules give rise to an interesting class of examples of adjoint pairs of contravariant functors. If $_RU_S$ is a partial cotilting bimodule, the functors in the adjoint pair $(\Hom_R(-, U), \Hom_S(-, U))$ have finite cohomological dimension; thus the derived functors $\R^b\Hom_R(-,U)$ and $\R^b\Hom_S(-, U)$ form a right adjoint pair in $\D^b(R\lMod)$ and $\D^b(\rMod S)$.  If $P$ is a projective module in $R\lMod$, $\Hom_R(P, U)$ is a direct summand of $U^{\alpha}_S$ for a suitable cardinal $\alpha$, and so, by condition (2) in Definition~\ref{def:cot}, $\Hom_R(P, U)$ is  $\Hom_S(-, U)$-acyclic. Thus $R\lMod$, and similarly $\rMod S$, have enough $\Hom_S(-, U)$-$\Hom_R(-, U)$-acyclic objects. Conversely, it is interesting to observe that, given a bimodule $_RU_S$, to assume both the finite cohomological dimension of $\Hom_R(-, U)$ and $\Hom_S(-, U))$, and the $\Hom_S(-, U)$-$\Hom_R(-, U)$-acyclicity of the projectives, implies that $_RU_S$ is a partial cotilting bimodule.

\section{Reflexive complexes}

Let us now investigate the relation between the $\D$-reflexivity of a complex in $\D^*(\mathcal A)$ and the $\D$-reflexivity of its terms or its cohomologies. This analysis will have an essential role in order to obtain our main results in the fourth and fifth sections. We always assume  that $(\Phi, \Psi)$ is an adjoint pair inducing the adjoint pair $(\R^*\Phi,\R^\dag\Psi)$.

\begin{thm}\label{prop:hom1}\label{thm:hom} 
Let $X$ be an object in $\D^*(\A)$.
\begin{enumerate}
\item If $X\in \D^b(\A)$ and any term of $X$ is $\D$-reflexive, then $X$ is $\D$-reflexive;
\item if $X\in \D^b(\A)$ and $H^i(X)$ is $\D$-reflexive for each $i\in\mathbb Z$, then $X$ is $\D$-reflexive.
\end{enumerate}
Assume $\R^\dag\Psi\R^*\Phi$ is way-out left (resp. right).
\begin{enumerate}
\item[(3)]  If  $X\in\D^-(\A)$ (resp. $\D^+(\A)$) and any term of $X$ is $\D$-reflexive, then $X$ is $\D$-reflexive;
\item[(4)] if $X\in\D^-(\A)$ (resp. $\D^+(\A)$) and $H^i(X)$ is $\D$-reflexive for any $i$, then $X$ is $\D$-reflexive.
\end{enumerate}
Assume $\R^\dag\Psi\R^*\Phi$ is way-out on both directions.
\begin{enumerate}
\item[(5)]  If any term of $X$ is $\D$-reflexive, then $X$ is $\D$-reflexive;
\item[(6)] if $H^i(X)$ is $\D$-reflexive for any $i$, then $X$ is $\D$-reflexive.
\end{enumerate}
\end{thm}
\begin{proof}
1. We can assume $X:= 0\to X_{-n}\to X_{-n+1}\to ...\to X_0\to 0$.
The thesis follows easily, by induction on the length $n$ of $X$, considering the triangles
 \[\tau_{\scriptscriptstyle> -1}(X)\to X \to \tau_{\scriptscriptstyle\leq -1}(X)\to
 \tau_{\scriptscriptstyle>-1}(X)[1].\]
 
2, 4, 6. The results follows applying \cite[Chp. I, Prop. 7.1]{H} to the morphism $\hat\eta\colon 1_{\D^*(\mathcal A)}\to \R^{\dag}\Psi\R^*\Phi$ and considering    the thick  subcategory of $\D$-reflexive objects of $\A$.

3. For short we denote by $\G$ the composition $\R^\dag\Psi\R^*\Phi$. We prove the result for $\G$ way-out left; the right case is analogous.
{\color{black} Following the proof of \cite[I.7.1]{H}, for each $j\in\mathbb Z$, it is possible to find a suitable $n\in\mathbb Z$ such that
\[H^j(\tau_{>n}X)\cong H^j(X)\text{ and }H^j(\G\tau_{>n}X)\cong H^j(\G X).\]
Then the conclusion follows since $H^j(\G\tau_{>n}X)\cong H^j(\tau_{>n}X)$ by part 1.
}

5. Let $X\in\D(\A)$; consider the triangle  \[\tau_{\scriptscriptstyle > 0} X \to X \to \tau_{\scriptscriptstyle \leq 0} X\to
 \tau_{\scriptscriptstyle > 0} X[1];\]
From 3 we know that $\tau_{\scriptscriptstyle \leq 0} X$ is $\D$-reflexive since $\R^\dag\Psi\R^*\Phi$ is way-out left and that $\tau_{\scriptscriptstyle > 0} X$ is $\D$-reflexive since $\R^\dag\Psi\R^*\Phi$ is way-out right. Thus we conclude that $X$ is $\D$-reflexive.
\end{proof}

The converse of the previous theorem is not in general true: in the following examples we show that there exist $\D$-reflexive complexes with not $\D$-reflexive terms  or not $\D$-reflexive cohomologies.

 \begin{EX}\label{ex:Gdim}
Let $\Lambda$ denote the $k$-algebra given by the quiver
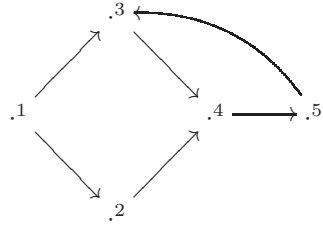
\[
\xymatrix{& {\cdot}^{3}\ar[dr]\\
{\cdot}^{1}\ar[ur] \ar[dr] && {\cdot}^{4}\ar[r]& {\cdot}^{5} \ar@/_13pt/[llu]\\
&{\cdot}^{2}\ar[ur]}
\]
with relations such that the left projective modules are $\begin{smallmatrix} & 1 & \\ 2 & & 3\\ & 4 &\end{smallmatrix}$, $\begin{smallmatrix} 2\\4\end{smallmatrix}$, $\begin{smallmatrix} 3\\4\\5\end{smallmatrix}$, $\begin{smallmatrix} 4\\5\end{smallmatrix}$ and $\begin{smallmatrix} 5\\3\end{smallmatrix}$. Consider the  module $_{\Lambda}U= \begin{smallmatrix} 5\end{smallmatrix}\oplus\begin{smallmatrix} 3\\4\\5 \end{smallmatrix}\oplus\begin{smallmatrix} & 1 & \\ 2 & & 3\\ & 4 &\end{smallmatrix} $ and let $S=\End_{\Lambda}(U)$. The algebra $S$ is given by the quiver $\cdot^{6}\to \cdot^{7}\to \cdot^{8}$ with right projectives $\begin{smallmatrix} 8\\7\end{smallmatrix}$, $\begin{smallmatrix} 7\\6\end{smallmatrix}$ and $\begin{smallmatrix} 6\end{smallmatrix}$, and $U_{S}=\begin{smallmatrix} 8\\7\end{smallmatrix}\oplus\begin{smallmatrix} 7\\6\end{smallmatrix}\oplus\begin{smallmatrix} 7\\6\end{smallmatrix}\oplus\begin{smallmatrix} 6\end{smallmatrix}\oplus\begin{smallmatrix} 6\end{smallmatrix}$.  Since $_{\Lambda}U_S$ is a partial cotilting bimodule,  $(\R^b\Hom_{\Lambda}(-,U), \R^b\Hom_{S}(-,U))$ is a  right adjunction and  the projective $\Lambda$-modules are $\Hom_{S}(-,U)$-$\Hom_{\Lambda}(-,U)$-acyclic objects.
Consider the complex with projective terms
\[P : \xymatrix{ 0 \ar[r] &  {\begin{smallmatrix} 4\\5\end{smallmatrix}} \ar[r] &  {\begin{smallmatrix} 3\\ 4\\5\end{smallmatrix}} \ar[r]  & {\begin{smallmatrix} 5\\3\end{smallmatrix}} \ar[r]  &  {\begin{smallmatrix} 3\\ 4\\5\end{smallmatrix}} \ar[r] & {\begin{smallmatrix} & 1 & \\ 2 & & 3\\ & 4 &\end{smallmatrix}} \ar[r] & 0}\]
and the obvious non-zero differentials.
It is easy to check that the morphism $\hat{\eta}_P$, given by the diagram
$$ \xymatrix{ 0 \ar[r] &  {\begin{smallmatrix} 4\\5\end{smallmatrix}} \ar[r] \ar[d]^{\eta_{P(4)}} &  {\begin{smallmatrix} 3\\ 4\\5\end{smallmatrix}} \ar[r] \ar[d]^{\eta_{P(3)}} & {\begin{smallmatrix} 5\\3\end{smallmatrix}} \ar[r] \ar[d]^{\eta_{P(5)}} &  {\begin{smallmatrix} 3\\ 4\\5\end{smallmatrix}} \ar[r] \ar[d]^{\eta_{P(3)}} & {\begin{smallmatrix} & 1 & \\ 2 & & 3\\ & 4 &\end{smallmatrix}} \ar[r] \ar[d]^{\eta_{P(1)}}& 0\\
0 \ar[r] &  {\begin{smallmatrix} 3\\ 4\\5\end{smallmatrix}} \ar[r]^{\cong} & {\begin{smallmatrix} 3\\4\\5\end{smallmatrix}} \ar[r]  &  {\begin{smallmatrix} 5\end{smallmatrix}} \ar[r] & {\begin{smallmatrix} 3\\ 4\\5\end{smallmatrix}} \ar[r]  &   {\begin{smallmatrix} & 1 & \\ 2 & & 3\\ & 4 &\end{smallmatrix}} \ar[r] &0}
$$
is a quasi-isomorphism. Nevertheless the terms $P(5)$ and $P(4)$ are not $\D$-reflexive.
\end{EX}

\begin{EX}\label{ex:GGdim}
Let $\Lambda$ be the $k$-algebra given by the quiver 

$$
\xymatrix{ {\cdot}1 \ar[r]
 & {\cdot}2  \ar[r] \ar@/_10pt/[l]
 & {\cdot}3 \ar[r] & {\cdot}4 \ar[r] & {\cdot}5  \ar@/_13pt/[ll] }
$$
with  relations such that the left projective modules are $\begin{smallmatrix} 1 \\ 2\\1\end{smallmatrix}$, $\begin{smallmatrix}  & 2 & \\1 & & 3 \\ &  & 4 \end{smallmatrix}$, $\begin{smallmatrix} 3 \\ 4\\5\end{smallmatrix}$, $\begin{smallmatrix} 4 \\ 5 \end{smallmatrix}$ and $\begin{smallmatrix} 5\\ 3 \end{smallmatrix}$. Let us consider the module $_{\Lambda}U=\begin{smallmatrix} 2\\1 \end{smallmatrix}\oplus \begin{smallmatrix} 1\\2\\1 \end{smallmatrix}\oplus \begin{smallmatrix} 5\\3 \end{smallmatrix}$ and let $S=\End_{\Lambda}(U)$. Then $S$ is given by the quiver
$$
\xymatrix{ {\cdot}6&  {\cdot}7 \ar[r]
 & {\cdot}8  \ar@/_10pt/[l]
 } 
$$ with relations such that the right projectives are   $\begin{smallmatrix} 8\\7\end{smallmatrix}$,  $\begin{smallmatrix} 7 \\ 8\\7\end{smallmatrix}$,  $\begin{smallmatrix} 6 \end{smallmatrix}$, and $U_{S}=\begin{smallmatrix} 8\\7 \end{smallmatrix}\oplus \begin{smallmatrix} 7\\8\\7 \end{smallmatrix}\oplus \begin{smallmatrix} 6 \end{smallmatrix}\oplus\begin{smallmatrix} 6 \end{smallmatrix}$. Since $_{\Lambda}U_S$ is a partial cotilting bimodule,  $(\R^b\Hom_{\Lambda}(-,U), \R^b\Hom_{S}(-,U))$ is a  right adjunction and  the projective $\Lambda$-modules are $\Hom_{S}(-,U)$-$\Hom_{\Lambda}(-,U)$-acyclic objects. Let us consider  the complex $X\in\D^b(\mathcal A)$ with projective terms
$$0\to \begin{smallmatrix} 1\\2 \\1 \end{smallmatrix}\stackrel{f}{\to} \begin{smallmatrix} 1\\2\\1 \end{smallmatrix}\stackrel{f}{\to} \begin{smallmatrix} 1\\2\\1 \end{smallmatrix}\to 0$$  
where $\Imm f=\soc P(1)$ and $\Ker f=\rad P(1)$. This complex is $\D$-reflexive: indeed $\R^b \Hom_S( \R^b \Hom_{\Lambda}(X, U), U)=\Hom_S(\Hom_{\Lambda}(X, U), U)$ and $\begin{smallmatrix} 1\\2\\1 \end{smallmatrix}=\Hom_{S}(\Hom_{\Lambda}(\begin{smallmatrix} 1\\2\\1 \end{smallmatrix}, U), U)$,  so   the evaluation map  $\hat{\eta}_{X}$ is trivially a quasi-isomorphism.
 Nevertheless the cohomology module $\Ker f/\imm f=S(2)$ is not $\D$-reflexive. In fact, let us consider a projective resolution of $S(2)$
 $$P: \ \ 0\to {\begin{smallmatrix}4\\5 \end{smallmatrix}}\to{\begin{smallmatrix} 3\\ 4\\5 \end{smallmatrix}}\to {\begin{smallmatrix} 5\\ 3 \end{smallmatrix}}\oplus{\begin{smallmatrix} 4\\5 \end{smallmatrix}}\to {\begin{smallmatrix} 3\\ 4\\5 \end{smallmatrix}}\oplus{\begin{smallmatrix} 3\\ 4\\5 \end{smallmatrix}}\to  \begin{smallmatrix}  & 2 & \\1 & & 3 \\ &  & 4 \end{smallmatrix}\oplus{\begin{smallmatrix} 5\\3 \end{smallmatrix}} \to {\begin{smallmatrix} 1\\ 2\\1 \end{smallmatrix}}\oplus{\begin{smallmatrix} 3\\ 4\\5 \end{smallmatrix}}\to \begin{smallmatrix}  & 2 & \\1 & & 3 \\ &  & 4 \end{smallmatrix}\to 0.$$  
 An easy computation shows that  $\R^b\Hom_S(\R^b\Hom_{\Lambda}( S(2), U), U)=\Hom_S(\Hom_{\Lambda}( P, U),U)$ is the complex
$$ \ 0\to 0 \to {\begin{smallmatrix} 5\\3 \end{smallmatrix} }\stackrel{\cong}{\to} {\begin{smallmatrix} 5\\3 \end{smallmatrix} }\to {\begin{smallmatrix} 5\\3 \end{smallmatrix} }\oplus{\begin{smallmatrix} 5\\3 \end{smallmatrix} }\to  {\begin{smallmatrix} 2\\ 1 \end{smallmatrix} }\oplus{\begin{smallmatrix} 5\\3 \end{smallmatrix} }\to  {\begin{smallmatrix} 1\\ 2\\ 1 \end{smallmatrix} } \oplus {\begin{smallmatrix} 5\\3 \end{smallmatrix} } \to  {\begin{smallmatrix} 2\\ 1 \end{smallmatrix} }\to 0$$
which has non zero cohomologies in degrees $0$ and ${-3}$. So $\hat{\eta}_{S(2)}$ is not a quasi-isomorphism and the module $S(2)$ is not $\D$-reflexive.\end{EX}

Given a finitely generated cotilting module of injective dimension $\leq 1$ over an Artin algebra, in \cite{CbCpF} is proved, using our terminology (see the forthcoming Theorem~\ref{cor:BB}) that the class of $\D$-reflexive modules coincides with the class of finitely generated ones. This can be generalized to cotilting modules of arbitrary finite injective dimension; in particular we obtain, in this setting, a converse of Theorem~\ref{thm:hom}. 

\begin{thm}\label{thm:artin2}
Let $\Lambda$ be an Artin algebra, $_{\Lambda}U$ a finitely generated cotilting module and $S=\End {}_{\Lambda}U$. Consider the adjoint pair $(\R^b\Hom_{\Lambda}(-, U), \R^b\Hom_S(-,U))$.
 \begin{enumerate}
 \item  A  complex $X\in\D^b(\Lambda\lMod)$ is  $\D$-reflexive if and only if the cohomologies $\Hm^{i}(X)$, $i\in\mathbb Z$, are finitely generated.
\item The subcategory of  bounded  $\D$-reflexive complexes is equivalent to $\D^b(\Lambda\lmod)$. 
\end{enumerate}
In particular a complex $X\in\D^b(\Lambda\lMod)$ is  $\D$-reflexive if and only if the cohomologies $\Hm^{i} (X)$, $i\in\mathbb Z$, are $\D$-reflexive.
\end{thm}
\begin{proof}
We recall that the assumptions imply that $_{\Lambda}U_S$ is a faithfully balanced cotilting bimodule (\cite[Theorem~1.5]{M}).

1. Let ${}^*=\Hom(-, W)$ be the usual duality between $\rmod\Lambda$ and $\Lambda\lmod$, where $W$ is the minimal injective cogenerator. Then $_{\Lambda}U^*=V_{\Lambda}$ is a finitely generated tilting module (see \cite{M}). Recall that a $\Lambda$-module is reflexive  with respect to the adjoint pair $({}^*, {}^*)$ if and only if it is finitely generated, and that the adjoint pair $(-\otimes^{\Le}_{S} V, \R\Hom_{\Lambda}(V, -))$ defines an equivalence between $\D^b(\rMod\Lambda)$ and $\D^b(\rMod S)$ (see \cite{CPS, Hp}). Let now  $X\in\D^b(\Lambda\lMod)$ be a  $\D$-reflexive complex  and let $P\in\K^-(\Lambda\lMod)$  be a complex of projective modules quasi-isomorphic to $X$. Then $P$ is quasi-isomorphic to $\Hom_{S}(\Hom_{\Lambda}(P, U), U)$. Using the standard adjunction formulas, since $U=V^*$,  we get that  
\[\Hom_{S}(\Hom_{\Lambda}(  P, U), U)=\Hom_{S}(\Hom_{\Lambda}( P, V^*), U)\cong
 \Hom_{S}(\Hom_{\Lambda}( V, P^* ), U)=\]
 \[= \Hom_{S}(\Hom_{\Lambda}( V, P^* ), V^*)\cong
  \Hom_{S}(\Hom_{\Lambda}( V, P^* )\otimes_S V, W)=\]
  \[=\Hom_{\Lambda}(\Hom_{\Lambda}(V, P^*)\otimes_{S} V,W).\] Moreover $\Hom_{\Lambda}(V, P^*)=\R\Hom_{\Lambda}(V, P^*)$ and, since   $\Hom_{\Lambda}(V, I)$ is $(-\otimes_{S}V)$-acyclic for any injective $\Lambda$-module $I$ \cite[Lemma 1.7]{M}, we obtain that 
\[\Hom_{\Lambda}(V, P^*)\otimes_{S} V=\R\Hom_{\Lambda}(V, P^*)\otimes^{\Le}_S V\cong P^*.\]
Hence $P$ is quasi-isomorphic to  $P^{**}$. For ${}^*$ is an exact functor, we conclude that $\Hm^{i} (P)$ is isomorphic to $ \Hm^{i} (P)^{**}$ for any $i$. Thus all the cohomologies of $X$ are finitely generated. 
Conversely, if all the cohomologies of $X$ are finitely generated, they are $\D$-reflexive: indeed all finitely generated projective $\Lambda$-modules are reflexive with respect to the adjoint pair $(\Hom_R(-, U), \Hom_S(-, U))$. Then we conclude by Proposition~\ref{prop:refl} and Theorem~\ref{thm:hom}.

2. It is well known that the subcategory of complexes in $\D^b(\Lambda\lMod)$ with finitely generated cohomologies is equivalent to $\D^b(\Lambda\lmod)$ (see \cite[Proposition~I.4.8]{H}).
\end{proof}

Limiting strongly the way-out dimensions, it is possible to prove that a complex is $\D$-reflexive if and only if its cohomologies are $\D$-reflexive in a more general setting.

\begin{prop}\label{thm:atmostone}
Let $X$ be an object of $\D^*(\A)$.
Suppose  the functor $\R^{\dag}\Psi\R^*\Phi$ to be way-out left of upper dimension $\leq 0$ and way-out right of lower dimension $\geq -1$. 
Then $X$  is $\D$-reflexive  if and only if its cohomologies are $\D$-reflexive.
\end{prop}
\begin{proof}For short, let us denote by $\G$ the composition $\R^{\dag}\Psi\R^*\Phi$.
 First let us suppose  $X\in \D^-(\mathcal A)$ to be a $\D$-reflexive complex. We can assume $X\in\D_{\leq 0}(\A)$  is of the form 
 \[X: \quad \dots {\to}  X_{-1}{\to}X_{0}\to 0.\]
 Let us first prove that  $H^{0}(X)$ is a $\D$-reflexive object.  Consider the triangle
 $$(*) \quad \sigma_{\scriptscriptstyle \leq -1 } X{\to} X{\to} \sigma_{\scriptscriptstyle > -1}X\to \sigma_{\scriptscriptstyle \leq -1} X[1] .$$ 
 The complex $\sigma_{\scriptscriptstyle > -1}X$ is quasi isomorphic to the stalk complex $H^{0}(X)$, and the complex  $\sigma_{\scriptscriptstyle \leq -1} X$ has zero cohomologies in degrees greater than $-1$.
Since $\dim^+\G\leq 0$, we have $H^i (\G X)=0$, $H^i(\G  \sigma_{\scriptscriptstyle > -1} X)=0$ and $H^i (\G \sigma_{\scriptscriptstyle \leq -1} X[-1])=H^{i-1} (\G \sigma_{\scriptscriptstyle \leq -1} X)=0$ for $i>0$.
Applying to the triangle $(*)$  first $\G$ and then the cohomology functor, we get the commutative diagram with exact rows
\small{$$ \xymatrix@-1.0pc{ \Hm^{-1}(\sigma_{\scriptscriptstyle \leq -1} X )\ar[r]^{\cong} \ar[d]  &  \Hm^{-1}(X) \ar[r] \ar[d]_{\cong}  & 0 \ar[r] \ar[d] &  0 \ar[r] \ar[d] & \Hm^{0}(X) \ar[r] \ar[d]_{\cong} & \Hm^{0}(\sigma_{\scriptscriptstyle > -1}X ) \ar[r] \ar[d] & 0\\
 \Hm^{-1}( \G \sigma_{\scriptscriptstyle \leq -1} X ) \ar[r] & \Hm^{-1}(\G  X) \ar[r] & \Hm^{-1} (\G  \sigma_{\scriptscriptstyle > -1}X) \ar[r] & 0 \ar[r] &   \Hm^{0}(\G X) \ar[r] & \Hm^{0}(\G  \sigma_{\scriptscriptstyle > -1}X ) \ar[r] & 0}$$}
\normalsize\noindent
Thus we deduce that $\Hm^{-1} (\G  \sigma_{\scriptscriptstyle > -1}X)=0$.   Since $\dim^-\G\geq -1$, we have $\Hm^{i} (\G \sigma_{\scriptscriptstyle > -1}X)=0$ for $i\leq -2$. Hence $\sigma_{\scriptscriptstyle > -1}X\cong \Hm^0(X)$ is $\D$-reflexive. Then, from the  triangle $(*)$  we deduce that  the complex $\sigma_{\scriptscriptstyle \leq -1 }X$ is $\D$-reflexive. Repeating the same argument for $\sigma_{\scriptscriptstyle \leq -1 } X[-1]$, we get that $H^{-1}(\sigma_{\scriptscriptstyle \leq -1 } X)\cong H^{-1} (X)$ is $\D$-reflexive. Continuing in such a way,  we conclude that $\Hm^i(X)$ is a $\D$-reflexive object for any $i\leq 0$. 

Suppose now $X$ to be a $\D$-reflexive complex in $\D(\A)$.
Consider the triangle
 $$\sigma_{\scriptscriptstyle \leq 0 } X{\to} X{\to} \sigma_{\scriptscriptstyle > 0}X\to \sigma_{\scriptscriptstyle \leq 0} X[1] $$ 
For the way-out dimensions of $\G$, we have   $H^i (\G \sigma_{\scriptscriptstyle >0} X )=0$ for $i<0$ and $H^i (\G \sigma_{\scriptscriptstyle \leq 0} X )=0$ for $i>0$. So  we get the commutative exact diagram
\footnotesize{$$ \xymatrix@-1.2pc{ \Hm^{-1}(X) \ar[r] \ar[d]_{\cong}  & 0 \ar[r] \ar[d] &  H^0(\sigma_{\scriptscriptstyle \leq 0} X ) \ar[r] \ar[d] & \Hm^{0}(X) \ar[r] \ar[d]_{\cong} & 0\ar[r] \ar[d] & 0 \ar[r] \ar[d] & H^1(X)\ar[r] \ar[d]^{\cong} & \Hm^{1}(\sigma_{\scriptscriptstyle > 0}X ) \ar[r] \ar[d] & 0\\
 \Hm^{-1}(\G  X) \ar[r] & 0\ar[r] & \Hm^{0} (\G  \sigma_{\scriptscriptstyle \leq 0}X) \ar[r] & H^0(\G X) \ar[r] & \Hm^{0} (\G  \sigma_{\scriptscriptstyle > 0}X) \ar[r]  &  0 \ar[r] &   \Hm^{1}(\G X) \ar[r] & \Hm^{1}(\G  \sigma_{\scriptscriptstyle > 0}X ) \ar[r] & 0}$$}
\normalsize
from which we conclude that $\sigma_{\scriptscriptstyle \leq 0 } X$ and $\sigma_{\scriptscriptstyle > 0}X$ are $\D$-reflexive complexes. Since the complex $\sigma_{\scriptscriptstyle \leq 0 } X$ belongs to $\D^-(\A)$, for what we have already proved we get that  $H^i(X)$ is $\D$-reflexive for any $i\leq 0$.
Similarly, considering the truncation in degree $i>0$ and the triangle
$$\sigma_{\scriptscriptstyle \leq i } X{\to} X{\to} \sigma_{\scriptscriptstyle > i}X\to \sigma_{\scriptscriptstyle \leq i} X[1] ,$$
we conclude  that $H^i(X)$ is $\D$-reflexive for any index $i$.
\end{proof}

\begin{cor}\label{cor:atmostone}
Assume  $\Phi$ has cohomological dimension at most one and $\A$ has enough $\Psi$-$\Phi$-acyclic objects.  If  $\mathcal B$ has enough $\Psi$-acyclics,   then  
 a complex $X\in \D(\mathcal A)$  is $\D$-reflexive if and only if its cohomologies $\Hm^i(X)$ are $\D$-reflexive.
\end{cor}
\begin{proof}
By Proposition~\ref{rem:dim}, $\R^{\dag}\Psi\R^*\Phi$ is way-out left of upper dimension  $\leq 0$ and way-out right of lower dimension $\geq -1$. So we can apply Proposition~\ref{thm:atmostone}.
\end{proof}

\begin{EX}
As in Example~\ref{ex:fasci}, let $(X, \mathcal O_X)$ be a locally noetherian scheme such that every coherent sheaf on $X$ is a quotient of a locally free sheaf  of finite rank.   Assume the structure sheaf   $\mathcal O_X$ has injective dimension one. Consider the abelian category $\mathfrak{Mod} X$ of sheaves of $\mathcal O_X$-modules and the thick subcategory $\mathfrak{Coh}X$ of coherent sheaves. Then  $(\R\mathcal{H}om(-, \mathcal O_X), \R\mathcal{H}om(-, \mathcal O_X))$ is a right adjoint pair in 
$\D({\mathfrak{Coh}X})$ which satisfies the assumptions of the previous corollary. Indeed, let $\mathcal L$ be the class of locally free sheaves of finite rank. Any object in $\mathcal L$ is $\mathcal{H}om(-, \mathcal O_X)$-acyclic and any $\mathcal F\in \mathfrak{Coh}X$ is image of a locally free sheaf of finite rank. Moreover, for any $\mathcal G$ locally free of finite rank, $\mathcal{H}om(\mathcal G, \mathcal O_X)$ is locally a finite direct sum of copies of $\mathcal O_X$ and so it is $\mathcal{H}om(-, \mathcal O_X)$-acyclic. Thus $\mathcal L$ satisfies the assumption of Definition~\ref{de:aciclici}.  Finally $\R\mathcal{H}om(-, \mathcal O_X)$ has cohomological dimension one.

Applying Proposition~\ref{prop:refl} we get that  any locally free sheaf of finite rank is $\D$-reflexive. Thus, considering   locally free resolutions, by Theorem~\ref{thm:hom} we  obtain that  any coherent sheaf is $\D$-reflexive and so   any complex in $\D({\mathfrak{Coh}X})$ is $\D$-reflexive. 
Note that Corollary~\ref{cor:atmostone} is trivially verified: indeed, if $\mathcal Y$ is a $\D$-reflexive complex in $\D({\mathfrak{Coh}X})$, then its cohomologies are $\D$-reflexive objects, being coherent sheaves (cfr. \cite[Prop. V. 2.1]{H}).
\end{EX}

The following technical result will be useful in the fifth section; its proof follows the same arguments used proving Proposition~\ref{thm:atmostone}.
{\color{black}
\begin{lemma}\label{prop:lastt}
Let $X$ be a $\D$-reflexive object in $\D^{*}_{\leq n}(\A)$, $n\in\mathbb Z$.
Suppose  the functor $\R^{\dag}\Psi\R^{*}\Phi$ to be way-out left of upper dimension $\leq 0$ and 
that $\R^*\Phi\Hm^j(X)$ is a stalk complex for each $j \in\mathbb Z$. For each $j\in\mathbb Z$, let $\rho(j)$ be an integer  such that
\[ H^{i}(\R^*\Phi\Hm^j(X))=0 \text{ for each $i\not=\rho(j)$.}\]
Then the cohomologies of $X$ are $\D$-reflexive if and only if 
\[H^i(\R^{\dag}\Psi H^{\rho(j)}(\R^*\Phi\Hm^j(X)))=0 \text{ for each }i\not=\rho(j), \rho(j)-1.\]  
\end{lemma}
\begin{proof}
 Consider the triangle
 \[
 (*) \quad \sigma_{\scriptscriptstyle \leq n-1 } X{\to} X{\to} \sigma_{\scriptscriptstyle > n-1}X\to \sigma_{\scriptscriptstyle \leq n-1} X[1] .
 \]
 The complex $\sigma_{\scriptscriptstyle > n-1}X$ is quasi isomorphic to the stalk complex $H^{n}(X)[-n]$, and the complex  $\sigma_{\scriptscriptstyle \leq n-1} X$ has zero cohomologies in degrees greater than $n-1$. By hypothesis, $H^i(\R^{*}\Phi H^n(X))=0$ for each $i\not=\rho(n)$; therefore $\R^{*}\Phi (H^{n}(X)[-\rho(n)])=\R^{*}\Phi H^{n}(X)[\rho(n)]$ is quasi isomorphic to the stalk complex $H^{\rho(n)}(\R^{*}\Phi H^n(X))$. Let us denote by $\G$ the composition $\R^{\dag}\Psi\R^{*}\Phi$; then we have
 \[\G(\sigma_{\scriptscriptstyle > n-1}X)=\G(H^{n}(X)[-n])=\R^{\dag}\Psi(\R^{*}\Phi H^{n}(X)[n])=\]
 \[=
  \R^{\dag}\Psi (H^{\rho(n)} (\R^{*}\Phi H^n(X))[n-\rho(n)])= \R^{\dag}\Psi H^{\rho(n)} (\R^{*}\Phi H^n(X))[\rho(n)-n].\]
 Since $\dim^+\G\leq 0$, we have $H^i (\G X)=0$, $H^i(\G  \sigma_{\scriptscriptstyle > n-1} X)=0$ and $H^i (\G \sigma_{\scriptscriptstyle \leq n-1} X[-1])=H^{i-1} (\G \sigma_{\scriptscriptstyle \leq n-1} X)=0$ for $i>n$. Applying to the triangle $(*)$  first $\G$ and then the cohomology functor, we get the commutative diagram with exact rows
\footnotesize{$$ \xymatrix@-1.0pc{ \Hm^{n-1}(\sigma_{\scriptscriptstyle \leq n-1} X )\ar[r]^-{\cong} \ar[d]  &  \Hm^{n-1}(X) \ar[r] \ar[d]_{\cong}  & 0 \ar[r] \ar[d] &  0 \ar[r] \ar[d] & \Hm^{n}(X) \ar[r] \ar[d]_{\cong} & \Hm^{n}(\sigma_{\scriptscriptstyle > n-1}X ) \ar[r] \ar[d] & 0\\
 \Hm^{n-1}( \G \sigma_{\scriptscriptstyle \leq n-1} X ) \ar[r] & \Hm^{n-1}(\G  X) \ar[r] & \Hm^{n-1} (\G  \sigma_{\scriptscriptstyle > n-1}X) \ar[r] & 0 \ar[r] &   \Hm^{n}(\G X) \ar[r] & \Hm^{n}(\G  \sigma_{\scriptscriptstyle > n-1}X ) \ar[r] & 0}$$}
\normalsize\noindent
Thus we deduce that $\Hm^{n-1} (\G  \sigma_{\scriptscriptstyle > n-1}X)=0$ and that $\Hm^{n} (\G  \sigma_{\scriptscriptstyle > n-1}X)\cong \Hm^n(\sigma_{\scriptscriptstyle > n-1}X)$. Since $\Hm^{n-i}(\sigma_{\scriptscriptstyle > n-1}X)=0$ for each $i>0$, the complex $\sigma_{\scriptscriptstyle > n-1}X$ is $\D$-reflexive, and hence the $n$-th cohomology of $X$ is $\D$-reflexive, if and only if for each $i>1$ we have
\[0=\Hm^{n-i} (\G  \sigma_{\scriptscriptstyle > n-1}X)=
\Hm^{n-i}(\G(H^{n}(X)[-n]))=\]
\[=\Hm^{n-i}(\R^{\dag}\Psi H^{\rho(n)}(\R^{*}\Phi H^n(X))[\rho(n)-n])=\Hm^{\rho(n)-i}\R^{\dag}\Psi H^{\rho(n)}(\R^{*}\Phi H^n(X)).
\]
Next, from the triangle $(*)$, $\sigma_{\scriptscriptstyle > n-1}X$ is $\D$-reflexive if and only if $\sigma_{\scriptscriptstyle \leq n-1}X$ is $\D$-reflexive. Applying the same argument to $\sigma_{\scriptscriptstyle \leq n-1}X$, we prove that $\Hm^{n-1} (X)$ is $\D$-reflexive if and only if the cohomologies $H^i(\R^{\dag}\Psi H^{\rho(n-1)}(\R^{*}\Phi H^{n-1}(X)))=0$ for each $i\not=\rho(n-1), \rho(n-1)-1$. Iterating this procedure, we conclude.
 \end{proof}

 \begin{rem}
Observe that if the functor $\Phi$ has cohomological dimension $\leq 1$, {\color{black} and there are enough $\Psi$-acyclic objects, 
under the hypotheses of Lemma~\ref{prop:lastt} the condition $H^i(\R^{\dag}\Psi H^{\rho(j)}(\R^{*}\Phi \Hm^j(X)))=0$ for each $i\not=\rho(j),\rho(j)-1$ is always satisfied (compare with Corollary~\ref{cor:atmostone}). The key point is that
\[R\Psi H^{\rho(j)}R\Phi H^jX\cong (R\Psi R\Phi H^jX)[-\rho(j)].
\]
We have $H^i R\Psi H^{\rho(j)}R\Phi H^jX=0$ for $i<0$ (because there are enough $\Psi$-acyclic objects) and
$H^i (R\Psi R\Phi H^jX)[-\rho(j)] =H^{i-\rho(j)} (R\Psi R\Phi H^jX)=0$ for $i>\rho(j)$ (because $R\Psi R\Phi$ is way out left of upper dimension $\leq 0$). Since $\Phi$ has cohomological dimension $\leq 1$, it is $|\rho(j)|\leq 1$ and so we conclude. }
\end{rem}
}
Let us see now the connection between the cohomological dimension of $\Phi$ and the closure of the class of $\D$-reflexive objects under kernels and cokernels.

\begin{thm}\label{prop:abelianita}
Assume that $\A$ has enough $\Psi$-$\Phi$-acyclic objects and $\B$ has enough $\Psi$-acyclic objects.
If $\Phi$ has cohomological dimension $\leq n$, then, in any exact sequence 
\[M_1\stackrel{f_1}{\to} M_2\stackrel{f_2}{\to} ...\stackrel{f_{n}}{\to} M_{n+1}\]
of $\D$-reflexive objects of $\mathcal A$, the kernels and the cokernels of the morphisms $f_i$, $i=1,...,n$, are $\D$-reflexive.

In particular, if $\Phi$ has cohomological dimension at most one, the class of $\D$-reflexive objects in $\mathcal A$ is an exact abelian subcategory of $\mathcal A$.
\end{thm}
\begin{proof}
First observe that by Proposition~\ref{rem:dim}, for any object $A$ in $\mathcal A$, the object $\R^{\dag}\Psi\R^*\Phi(A)$ belongs to
$\mathcal D_{\geq -n}(\mathcal A)\cap \mathcal D_{\leq 0}(\mathcal A)$.
Denoted by $K_i$ the kernel of the morphism $f_i$, $i=1,...,n$, by $K_{n+1}$ the image of $f_{n}$, and by $K_{n+2}$ the cokernel of $f_{n}$,  let us consider the following triangles in $\D^b(\mathcal A)$:
\[K_i\to M_i\to K_{i+1}\to K_i[1],\quad i=1,...,n+1.\]
Consider the maps $H^0(\hat\eta_{M_i}):M_i\to H^0(\R^\dag\Psi\R^*\Phi M_i)$, and $H^0(\hat\eta_{K_j}):K_j\to H^0(\R^\dag\Psi\R^*\Phi K_j)$, $1\leq i\leq n$, $1\leq j\leq n+2$. Since $M_i$ are $\D$-reflexive objects, clearly $H^0(\hat\eta_{M_i})$ are isomorphisms.
We will prove that $H^j(\R^\dag\Psi\R^*\Phi K_i)=0$ for each $j\not=0$ and each $1\leq i\leq n+2$ and that all $H^0(\hat\eta_{K_j})$ are isomorphisms.

Because of the way-out dimension of $\R^\dag\Psi\R^*\Phi$, $H^j(\R^\dag\Psi\R^*\Phi K_i)=0$ for each $j>0$.
Applying the cohomology functor we get the long exact sequences
\[0\to H^{-n}(\R^\dag\Psi\R^*\Phi K_i)\to 0\to H^{-n}(\R^\dag\Psi\R^*\Phi K_{i+1})\to
H^{-n+1}(\R^\dag\Psi\R^*\Phi K_i)\to 0\to ...\]
\[...\to 0\to  H^{-2}(\R^\dag\Psi\R^*\Phi K_{i+1})\to H^{-1}(\R^\dag\Psi\R^*\Phi K_{i})\to 0\to\]
\[\to H^{-1}(\R^\dag\Psi\R^*\Phi K_{i+1})\to H^{0}(\R^\dag\Psi\R^*\Phi K_{i})\to H^{0}(\R^\dag\Psi\R^*\Phi M_{i})\to H^{0}(\R^\dag\Psi\R^*\Phi K_{i+1})\to 0,\]
for $i=1,...,n+1$.
In particular $H^{-n}(\R^\dag\Psi\R^*\Phi K_i)=0$ for $i=1,..., n+1$; therefore for $j=1,..., n$, since $n-j+1<n-j+2\leq n+1$, we have
\[H^{-j}(\R^\dag\Psi\R^*\Phi K_1)\cong H^{-n}(\R^\dag\Psi\R^*\Phi K_{n-j+1})=0, \text{ and}\]
\[H^{-j}(\R^\dag\Psi\R^*\Phi K_2)\cong H^{-n}(\R^\dag\Psi\R^*\Phi K_{n-j+2})=0.\]
Working a little on diagrams
\[\xymatrix@-1pc{
&0\ar[r] &K_i\ar[r]\ar[d]^{H^0(\hat\eta_{K_i})} &M_i\ar[r]\ar[d]^{\cong} &K_{i+1}\ar[r]\ar[d]^{H^0(\hat\eta_{K_{i+1}})} &0\\
&...\ar[r] &H^{0}(\R^\dag\Psi\R^*\Phi K_i)\ar[r] &H^{0}(\R^\dag\Psi\R^*\Phi M_i)\ar[r] &H^{0}(\R^\dag\Psi\R^*\Phi K_{i+1})\ar[r] &0}
\]
with $i=1,2$,
we get that $H^0(\hat\eta_{K_1})$ and $H^0(\hat\eta_{K_2})$ are isomorphisms. Therefore $K_1$ and $K_2$ are $\D$-reflexive. Working with the triangles
\[K_i\to M_i\to K_{i+1}\to K_i[1],\quad i=2,...,n+1\]
we get that also $K_3$, ..., $K_{n+2}$ are $\D$-reflexive.
\end{proof}

\section{The 1-dimensional case}

In the previous section we have seen that more precise results are available when the involved functors have cohomological dimension at most one. This section is dedicated to study in detail this favorable case. Our aim is to characterize the $\D$-reflexive objects in the abelian categories $\A$ and $\B$, producing a general form of the {\color{black}Cotilting Theorem in the sense of Colby and Fuller (see \cite[Ch.~5]{CbF1}), a contravariant version of the celebrated Brenner and Butler Theorem \cite{BB}}.

We assume $\A$ has enough $\Psi$-$\Phi$-acyclic objects and $\B$ has enough $\Phi$-$\Psi$-acyclic objects, respectively, and $(\Phi,\Psi)$ is an adjoint pair of contravariant functors of cohomological dimension  at most one. In particular, under these assumptions,
\begin{itemize}
\item there exist the total derived functors $\R\Phi$ and $\R\Psi$, and they have both lower dimension $\geq 0$ and upper dimension $\leq 1$,
\item the composition $\R\Psi\R\Phi$ results to be way-out left of upper dimension $\leq 0$ and way-out right of lower dimension $\geq -1$ (Proposition~\ref{rem:dim}), and it is isomorphic  to $\R(\Psi\Phi)$ (\cite[Proposition 5.4]{H}),
\item the families of $\Phi$-acyclic and $\Psi$-acyclic objects are closed under submodules (Proposition~\ref{lemma:aciclici}),
\item a complex is $\D$-reflexive if and only if its cohomologies are $\D$-reflexive (Corollary~\ref{cor:atmostone}),
\item the classes of $\D$-reflexive objects in $\A$ and $\B$ are exact abelian subcategories of $\A$ and $\B$ (Theorem~\ref{prop:abelianita}).
\end{itemize}

\noindent In this setting we deal with the unbounded derived categories $\D(\A)$ and $\D(\B)$ and the total derived functors $\R\Phi$ and $\R\Psi$: for any complex $X$ in $\D(\A)$ (resp. $\D(\B)$), we denote by $R^i\Phi X$ (resp. $R^i\Psi X$), $i\in\mathbb Z$, the $i^{th}$-cohomology $H^i(\R\Phi X)$ (resp. $H^i(\R\Psi X)$). Observe that $R^0\Phi A=\Phi A$ and $R^0\Psi B=\Psi B$ for each $A$ in $\A$ and $B$ in $\B$.
\begin{lemma}\label{lemma:Psiaciclici}
Any object in $\Imm \Phi$ is $\Psi$-acyclic.
\end{lemma}
\begin{proof}
 Let $A$ be an object in $\mathcal A$. Consider an epimorphism $L\to A\to 0$ where $L$ is a $\Psi$-$\Phi$-acyclic object.
Applying $\Phi$ we get the monomorphism $0\to \Phi A\to \Phi L$. Since $\Phi L$ is $\Psi$-acyclic, and the family of $\Psi$-acyclic objects is closed under submodules, we conclude that $\Phi A$ is $\Psi$-acyclic.
\end{proof}

Observe that by the previous lemma any $\Phi$-acyclic object is also $\Psi$-$\Phi$-acyclic.

\begin{prop}\label{cor:driflessivi}

An object $A\in\mathcal A$ is $\D$-reflexive if and only if $\Psi R^1\Phi A=0$ and the map $H^0(\hat\eta_A):A\to R^0(\Psi\Phi)A$ is an isomorphism.
\end{prop}
\begin{proof}  
Since $\R\Psi\R\Phi$ is way-out  of upper dimension $\leq 0$ and lower dimension $\geq -1$, the object $A$ is $\D$-reflexive if and only if $H^0(\hat\eta_A)$ and $H^{-1}(\hat\eta_A)$ are isomorphisms, the latter being equivalent to $H^{-1}(\R(\Psi\Phi)A)=0$.
Let us consider the triangle
\[
\sigma_{\scriptscriptstyle \leq 0}\R\Phi A\to \R\Phi A\to\sigma_{\scriptscriptstyle > 0}\R\Phi A \to \sigma_{\scriptscriptstyle \leq 0}\R\Phi A[1];
\]
taking in account the way-out dimensions of $\R\Phi$, this triangle is isomorphic to
\[\Phi A\to \R\Phi A\to  R^1\Phi A[-1] \to \Phi A [1]\]
Applying $\R\Psi$ we get, using Lemma~\ref{lemma:Psiaciclici}, the triangle
\[\R\Psi( R^1\Phi A [-1])\to\R\Psi\R\Phi A=\R(\Psi\Phi)A\to \Psi \Phi A \to \R\Psi  R^1\Phi A[2] .
\]
Considering the associated cohomology sequence, we get the exact sequence
\[0\to H^{-1}(\R\Psi( R^1\Phi A[-1]))\to H^{-1}(\R(\Psi\Phi)A)\to 0;\]
Then we conclude since
\[H^{-1}(\R\Psi( R^1\Phi A[-1]))
=H^0(\R\Psi  R^1\Phi A)=\Psi R^1\Phi A.\]
\end{proof}

\begin{thm}\label{cor:BB}
An object $A$ in $\A$ is $\D$-reflexive if and only if
\begin{enumerate}
\item $\Phi(A)$ and $ R^1\Phi(A)$ are $\D$-reflexive;
\item $ R^i\Psi R^j\Phi(A)=0$ if $i\not=j$;
\item there exists a natural map $\gamma_A$ and an exact sequence
$$0\to  R^1\Psi R^{1}\Phi(A)\stackrel{\gamma_A}{\to} A\stackrel{\eta_A}{\to} \Psi\Phi(A)\to 0$$
\end{enumerate}
In such a case, when denoting by $\pi_{\R\Phi A}$ the natural map $\R\Phi A\to
\sigma_{\scriptscriptstyle > 0}\R\Phi A$, we have $\gamma_A=H^0(\hat\eta_A)^{-1}\circ R^0\Psi(\pi_{\R\Phi A})$.
\end{thm}
\begin{proof}
Assume $A$ is $\D$-reflexive. Since $\R\Phi A$ is $\D$-reflexive, from Corollary~\ref{cor:atmostone} it follows that its cohomologies $\Phi(A)$ and $ R^1\Phi A$ are $\D$-reflexive. By Lemma~\ref{lemma:Psiaciclici}  $\Phi(A)$ is $\Psi$-acyclic; therefore
\[
R^1\Psi(R^0\Phi(A))=R^1\Psi(\Phi(A))=0.\] By Proposition~\ref{cor:driflessivi},  we know also  that $\Psi R^1\Phi(A)=0$, and so $R^0\Psi R^1\Phi(A)=0$.
To prove (3), let us consider the triangle
\[\sigma_{\scriptscriptstyle \leq 0}\R\Phi A\stackrel{\iota}{\to}\R\Phi A\stackrel{\pi_{\R\Phi A}}{\to}
\sigma_{\scriptscriptstyle > 0}\R\Phi A\to \sigma_{\scriptscriptstyle \leq 0}\R\Phi A[1];\]
taking in account the way-out dimensions of $\R\Phi$, this triangle is isomorphic to
\[
\Phi(A)\stackrel{\iota}{\to}\R\Phi A\stackrel{\pi_{\R\Phi A}}{\to}R^1\Phi(A)[-1]\to\Phi(A)[1].
\]
Applying $\R\Psi$, by Lemma~\ref{lemma:Psiaciclici} we get the triangle
\[
\R\Psi( R^1\Phi(A))[1]\stackrel{\R\Psi(\pi_{\R\Phi A})}{\to}\R\Psi\R\Phi A\stackrel{\R\Psi(\iota)}{\to}
\Psi\Phi(A)\to \R\Psi( R^1\Phi(A))[2].
\]
Considering the associated cohomology sequence, we get the natural short exact sequence
\[0\to R^1\Psi( R^1\Phi(A)) \stackrel{R^0\Psi(\pi_{\R\Phi A})}{\to}
H^0(\R\Psi\R\Phi A)=H^0(\R(\Psi\Phi)( A))\stackrel{R^0\Psi(\iota)}{\to}  \Psi(\Phi(A))\to 0.\]
Since $A$ is $\D$-reflexive, $H^0(\hat\eta_A):A\to H^0(\R(\Psi\Phi)( A))$ is an isomorphism. Denote by $\gamma_A$ the composition $H^0(\hat\eta_A)^{-1}\circ R^0\Psi(\pi_{\R\Phi A})$; we can apply Proposition~\ref{prop:legame} to get
$R^0\Psi(\iota)\circ H^0(\hat\eta_A)=\eta_A$
and hence the natural exact sequence
\[0\to  R^1\Psi R^1\Phi(A)\stackrel{\gamma_A}{\to} A \stackrel{\eta_A}{\to} \Psi\Phi(A)\to 0.\]
Conversely, assume conditions (1), (2) and (3) hold. Applying (1) to $\Phi(A)$ and $R^1\Phi(A)$, we get that $ R^1\Psi R^1\Phi(A)$ and $\Psi\Phi(A)$ are $\D$-reflexive. Therefore, by (3) also $A$ is $\D$-reflexive.
\end{proof}
{\color{black}The same result holds for any $\D$-reflexive object $B$ in $\B$, with the map $\theta_B:R^1\Phi R^{1}\Psi(B)\to B$, $\theta_B=H^0(\hat\xi_B)^{-1}\circ R^0\Phi(\pi_{\R\Psi B})$, which plays the role of the natural map $\gamma$.}

We are now ready to give a Cotilting Theorem in the sense of \cite[Ch.~5]{CbF1}, between the classes of $\D$-reflexive objects induced by the pair of adjoint functors $(\Phi,\Psi)$.

\begin{cor}\label{cor:CotiltingTh}
{\color{black}Consider the following subclasses of the abelian subcategories $\D_\A$ and $\D_\B$ of $\D$-reflexive objects in $\A$ and $\B$:
\[\mathcal T_{\A}=\Ker\Phi\cap\D_\A, \quad \mathcal F_{\A}=\Ker R^1\Phi\cap\D_\A\] 
\[\mathcal T_{\B}=\Ker\Psi\cap\D_\B, \quad \mathcal F_{\B}=\Ker R^1\Psi\cap\D_\A\]
Then the following conditions are satisfied:
\begin{enumerate}
\item $\Phi:\D_\A\to \F_\B$, $R^1\Phi:\D_\A\to \T_\B$, 
$\Psi:\D_\B\to \F_\A$, $R^1\Psi:\D_\B\to \T_\A$.
\item for each object $A$ in $\D_\A$ and $B$ in $\D_\B$ we have the following exact sequences of natural maps
\[0\to  R^1\Psi R^{1}\Phi(A)\stackrel{\gamma_A}{\to} A\stackrel{\eta_A}{\to} \Psi\Phi(A)\to 0\]
\[0\to  R^1\Phi R^{1}\Psi(B)\stackrel{\theta_B}{\to} B\stackrel{\xi_B}{\to} \Phi\Psi(B)\to 0\]
\item the restrictions
\[\Phi:\F_\A\dual \F_\B:\Psi\text{ and } R^1\Phi:\T_\A\dual \T_\B:R^1\Psi\]
define category equivalences.
\end{enumerate}
Moreover these are the largest possible classes where such a duality arises
}\end{cor}

Our starting point was that $(\Phi,\Psi)$ is an adjoint pair of functors between the abelian categories $\A$ and $\B$. Now we show that
$( R^1\Phi,  R^1\Psi)$ is an adjoint pair of functors between the abelian categories of $\D$-reflexive objects of $\A$ and $\B$.

\begin{thm}\label{cor:adj}
In the classes $\D_{\A}$ and $\D_{\B}$ of $\D$-reflexive objects of $\A$ and $\B$, the pair $( R^1\Phi,  R^1\Psi)$ is left adjoint with the natural maps $\gamma$ and $\theta$ as units.
\end{thm}
\begin{proof} 
In order to prove that  $( R^1\Phi,  R^1\Psi)$ is a left adjoint pair in the classes $\D_{\A}$ and $\D_{\B}$, it is enough to show that $\theta_{ R^1\Phi A}\circ  R^1\Phi(\gamma_A)=\id_{ R^1\Phi A}$ for any $A\in \D_{\A}$ and, analogously, $\gamma_{ R^1\Psi B}\circ  R^1\Psi(\theta_B)=\id_{ R^1\Psi B}$ for any $B\in \D_{\B}$. 

Note that, from the adjunction formula $\R\Phi(\hat\eta_A)\circ\hat\xi_{\R\Phi A}=\id_{\R\Phi A}$, we get that $R^1\Phi(\hat\eta_A)\circ H^1(\hat\xi_{\R\Phi A})=\id_{ R^1\Phi A}$. We will prove that $\theta_{ R^1\Phi A}=H^1(\hat\xi_{\R\Phi A})^{-1}$ and $ R^1\Phi (\gamma_A)=R^1\Phi(\hat\eta_A)^{-1}$.

First, let us consider the diagram
\[
\xymatrix{
{}\Phi A\ar[r]\ar[d]^{\hat\xi_{\Phi A}}&\R\Phi A\ar[r]^{\pi_{\R\Phi A}}\ar[d]^{\hat\xi_{\R\Phi A}}&R^1\Phi A[-1]\ar[r]\ar[d]^{\hat\xi_{R^1\Phi A[-1]}}&\Phi A[1]\ar[d]^{\hat\xi_{\Phi A[1]}}\\
{}\R\Phi\R\Psi({}\Phi A)\ar[r]&\R\Phi\R\Psi(\R\Phi A)\ar[r]&\R\Phi\R\Psi(R^1\Phi A[-1])\ar[r]&\R\Phi\R\Psi(\Phi A [1])}
\]
Applying the cohomology functor $H^1$ we get
\[
\xymatrix{
0\ar[r] &R^1\Phi A\ar[dd]^{H^1(\hat\xi_{\R\Phi A})}\ar@{=}[rrr]&&&R^1\Phi A\ar[dd]^{H^1(\hat\xi_{R^1\Phi A[-1]})}\ar[r]&0\\
(*)\\
0\ar[r] &H^1(\R\Phi\R\Psi(\R\Phi A))\ar[rrr]^-{H^1(\R\Phi\R\Psi(\pi_{\R\Phi A}))}&&&H^1(\R\Phi\R\Psi(R^1\Phi A[-1]))\ar[r]&0
}
\]
Let us prove that $H^1(\R\Phi\R\Psi(\pi_{\R\Phi A}))$ is the identity map. Consider a $\Phi$-acyclic resolution $P$ of $A$; then we have
\[
\xymatrix{
{}\R\Phi A:=\ar[d]_{\pi_{\R\Phi A}}&0\ar[r]&\Phi(P_0)\ar[d]^p\ar[r]^{\alpha}&\Phi(P_1)\ar[r]\ar@{=}[d]&...\\
\sigma_{\scriptscriptstyle > 0}\R\Phi A:=&0\ar[r]&\Phi(P_0)/\Ker\alpha\ar[r]&\Phi(P_1)\ar[r]&...}
\]
Since the terms in both the complexes are $\Phi$-$\Psi$-acyclic, we get
\[
\xymatrix{
\R\Phi\R\Psi\R\Phi A:=\ar[d]_{\R\Phi\R\Psi(\pi_{\R\Phi A})}&0\ar[r]&\Phi\Psi\Phi(P_0)\ar[d]^{\Phi\Psi(p)}\ar[r]^{\Phi\Psi(\alpha)}&\Phi\Psi\Phi(P_1)\ar[r]\ar@{=}[d]&...\\
\R\Phi\R\Psi(\sigma_{\scriptscriptstyle > 0}\R\Phi A)=&0\ar[r]&\Phi\Psi(\Phi(P_0)/\Ker\alpha)\ar[r]&\Phi\Psi\Phi(P_1)\ar[r]&...}
\]
Since the functor $\Phi\Psi$ is exact on the short exact sequence of $\Phi$-$\Psi$-acyclic objects
\[0\to \Ker \alpha\to \Phi(P_0) \stackrel p{\to} \Phi(P_0)/\Ker \alpha\to 0,\]
the map $\Psi\Phi(p)$ is surjective; it is now clear that 
\[H^1(\R\Phi\R\Psi(\pi_{\R\Phi A}))=1_{H^1(\R\Phi\R\Psi\R\Phi A)}.
\]
Since $\pi_{\R\Psi(R^1\Phi A)}:\R\Psi(R^1\Phi A)\to \sigma_{\scriptscriptstyle > 0}\R\Psi(R^1\Phi A)$ is the identity map, by
Theorem~\ref{cor:BB} and diagram $(*)$ we have
\[\theta_{R^1\Phi A}=H^0(\hat\xi_{R^1\Phi A})^{-1}=H^1(\hat\xi_{R^1\Phi A[-1]})^{-1}=
H^1(\hat\xi_{\R\Phi A})^{-1}.\]

Second, thinking at $\gamma_A:R^1\Psi R^1\Phi A\to A$ as a map between stalk complexes, let us consider the following commutative diagram (see Theorem~\ref{cor:BB})
\[
\xymatrix{
{}\R\Psi(\sigma_{\scriptscriptstyle > 0}\R\Phi A)\ar[rr]^-{\R\Psi(\pi_{\R\Phi A})}
\ar[d]^{\cong}_{qiso}
&&
\R\Psi\R\Phi A\ar[r]^-{{\hat\eta_A}^{-1}}&A
\\
\R\Psi(R^1\Phi A[-1])
\ar@{=}[d]\\
\R\Psi(R^1\Phi A)[1]
\ar[d]^{\cong}_{qiso}\\
R^1\Psi R^1\Phi A\ar[rrruuu]_{\gamma_A}
}
\]
Applying $\R\Phi$ we get the commutative diagram
\[
\xymatrix{
{}\R\Phi\R\Psi\R\Phi A\ar[rr]^-{\R\Phi\R\Psi(\pi_{\R\Phi A})}\ar[d]^{\R\Phi(\hat\eta_A)}&&
\R\Phi\R\Psi(\sigma_{\scriptscriptstyle > 0}\R\Phi A)\\
\R\Phi A\ar[rru]_{\R\Phi(\gamma_A)}
}
\]
Applying the cohomology functor $H^1$ we get
\[R^1\Phi(\gamma_A)\circ H^1(\R\Phi(\hat\eta_A))=H^1(\R\Phi\R\Psi(\pi_{\R\Phi A}))=1_{H^1(\R\Phi\R\Psi\R\Phi A)}.\]
\end{proof}

\begin{EX}
Let $R$ and $S$ be arbitrary associative rings. Consider a partial cotilting bimodule $_RU_S$  of injective dimension $\leq 1$. Then the adjunction $(\Hom_R(-,U), \Hom_S(-, U))$ satisfies  the assumptions of Theorem~\ref{cor:BB} and Theorem~\ref{cor:adj} . 
In particular the classes of $\D$-reflexive modules are abelian subcategories of $R\lMod$ and $\rMod S$ and the restriction of  the functors $\Ext^1(-,U)$ to these classes forms a left adjoint pair. Moreover,  {\color{black} following Corollary~\ref{cor:CotiltingTh},} the functors $\Ext^1(-,U)$ induce a duality between the subcategories of {\color{black}$\D$-reflexive modules in} $\Ker\Hom(-,U)$  and the functors $\Hom(-,U)$ between the subcategories of {\color{black}$\D$-reflexive modules in } $\Ker\Ext^1(-,U)$ (compare with \cite{Cb, Cb1, CbF, C, CF, Ma, T}).
\end{EX}

\begin{EX}
As in Example~\ref{ex:fasci}, let $(X, \mathcal O_X)$ be a locally noetherian scheme such that every coherent sheaf on $X$ is a quotient of a locally free sheaf.   Assume the structure sheaf   $\mathcal O_X$ has injective dimension one and consider the adjunction  $(\R\mathcal{H}om(-, \mathcal O_X), \R\mathcal{H}om(-, \mathcal O_X))$  in 
$\D({\mathfrak{Coh} X})$. As we have already seen, any coherent sheaf is  $\D$-reflexive. Since any  locally free sheaf of finite rank is $\mathcal{H}om(-, \mathcal O_X)$-$\mathcal{H}om(-, \mathcal O_X)$-acyclic, the assumption of Theorem~\ref{cor:BB} and Theorem~\ref{cor:adj} are satisfied. Denoted as $\mathcal{E}xt(-, \mathcal O_X)$ the first derived functor of $\mathcal{H}om(-, \mathcal O_X)$ (see \cite[Chp. III]{H^1}), we get that $(\mathcal{E}xt(-, \mathcal O_X), \mathcal{E}xt(-, \mathcal O_X))$ is a left adjoint pair in $\mathfrak{Coh}X$,   the functors $\mathcal{E}xt(-, \mathcal O_X)$ induce a duality between the coherent sheaves in $\Ker\mathcal{H}om(-, \mathcal O_X)$  and the functors $\mathcal{H}om(-, \mathcal O_X)$ between the coherent sheaves in $\Ker\mathcal{E}xt(-, \mathcal O_X)$.
\end{EX}

\section{The $n$-dimensional case}

In this section we recover a Cotilting Theorem in the case of functors of cohomological dimension greater than one. 
We assume $\A$ and $\B$ have enough projectives, $(\Phi,\Psi)$ is an adjoint pair of contravariant functors of cohomological dimension  at most $n$, $\Phi(P)$ is $\Psi$-acyclic for each projective $P$ in $\A$, and $\Psi(Q)$ is $\Phi$-acyclic for each projective $Q$ in $\B$. For instance this is the case when $\Phi$ and $\Psi$ are the contravariant $\Hom$-functors associated to a partial cotilting bimodule.

Let $P$ be a projective resolution of an object $A$ in $\A$. Denote by $Q_{**}$ a Cartan-Eilenberg resolution of $\Phi(P)$; applying $\Psi$ to the bicomplex $Q_{**}$, we get the bicomplex
\[
\xymatrix{
&...&...&...&...\\
...\ar[r] &\Psi Q_{2,-1}\ar[u]\ar[r]&\Psi Q_{1,-1}\ar[u]\ar[r]&\Psi Q_{0,-1}\ar[u]\ar[r]&0\\
...\ar[r]&\Psi Q_{2,0}\ar[u]\ar[r]&\Psi Q_{1,0}\ar[u]\ar[r]&\Psi Q_{0,0}\ar[u]\ar[r]& 0\\
&0\ar[u]&0\ar[u]&0\ar[u]
}
\]
To this bicomplex we associate two spectral sequences ${}_IE_2^{pq}$ and ${}_{II}E_2^{pq}$:
\[{}_IE_2^{pq}=H^p_h(H^q_v(\Psi Q_{**}))=\left\{\begin{array}{ll}
\label{}
    0   & \text{ if }q\not=0\\
  H^p_h(\Psi\Phi P)= R^p(\Psi\Phi)(A) &   \text{ if }q=0
\end{array}
\right.
\]
\[{}_{II}E_2^{pq}=H^p_v(H^q_h(\Psi Q_{**}))=H^p_v(\Psi(H^{-q}_h(Q_{**})))=H^p_v(\R\Psi(R^{-q}\Phi A))=
R^{p}\Psi R^{-q}\Phi(A).\]
Observe that ${}_IE_2^{pq}=0$ for either $p>0$ or $q\not=0$ and ${}_{II}E_2^{pq}=0$
for either $p<0$ or $q>0$.

Both these spectral sequences converge to the hypercohomology $\mathbb R^{q+p}\Psi(\R\Phi(A))$. The first spectral sequence ${}_IE_2^{pq}$ collapses to yield
\[\mathbb R^{n}\Psi(\R\Phi(A))=R^{n}(\Psi\Phi)(A),\]
which is zero for $n>0$.
The second spectral sequence ${}_{II}E_2^{pq}$ lies on the fourth quadrant:
\[
\xymatrix{
\ar@{--}[rrrr]&&&&\\
&R^0\Psi R^0\Phi(A)&R^1\Psi R^0\Phi(A)&R^2\Psi R^0\Phi(A)&...\\
&R^0\Psi R^{1}\Phi(A)&R^1\Psi R^{1}\Phi(A)&R^2\Psi R^{1}\Phi(A)&...\\
&R^0\Psi R^{2}\Phi(A)&R^1\Psi R^{2}\Phi(A)&R^2\Psi R^{2}\Phi(A)&...\\
\ar@{--}[uuuu]&...&...&...&
}
\]
together with maps
\[d_2^{pq}:{}_{II}E_2^{pq}=R^{p}\Psi R^{-q}\Phi(A)\to {}_{II}E_2^{p+2,q-1}=R^{p+2}\Psi R^{1-q}\Phi(A).\] 

Since the cohomological dimension of $\Phi$ and $\Psi$ is at most $n$, we have ${}_{II}E_{n+1}^{pq}={}_{II}E_{n+i}^{pq}={}_{II}E_{\infty}^{pq}$ for each $p$, $q$. For each $s\leq 0$, $R^{s}(\Psi\Phi)(A)$ has a finite filtration
\[0=F^{n+1+s}R^{s}(\Psi\Phi)(A)\subseteq F^{n+s}R^{s}(\Psi\Phi)(A)\subseteq ...\]
\[...\subseteq F^{1}R^{s}(\Psi\Phi)(A) \subseteq F^{0}R^{s}(\Psi\Phi)(A)=R^{s}(\Psi\Phi)(A)\]
with $[F^{i}R^{s}(\Psi\Phi)(A)]/[F^{i+1}R^{s}(\Psi\Phi)(A)]\cong {}_{II}E_{\infty}^{i,s-i}$.

If $\mathbf{n=1}$ we have ${}_{II}E_2^{pq}={}_{II}E_{\infty}^{pq}$ for each $p$ and $q$; therefore
\[
0=R^1(\Psi\Phi)(A)={}_{II}E_2^{10}=R^1\Psi R^0\Phi(A).
\]
If $A$ is $\D$-reflexive, $R^{-1}(\Psi\Phi)(A)=0$ and $R^0(\Psi\Phi)(A)\cong A$; hence using the edge homomorphisms, it is easy to get
\begin{enumerate}
\item ${}_{II}E_2^{0-1}=R^0\Psi R^{1}\Phi(A)=0$;
\item there exists the following short exact sequence with natural maps
\[0\to {}_{II}E_2^{1-1}=R^1\Psi R^{1}\Phi(A)\to A\to {}_{II}E_2^{00}=R^0\Psi R^0\Phi(A)\to 0.\]
\end{enumerate}
It is not hard now to recover Proposition~\ref{cor:driflessivi} and partially Theorem~\ref{cor:BB}, (3).

If $\mathbf{n=2}$, we have ${}_{II}E_2^{pq}={}_{II}E_{\infty}^{pq}$ for $(p,q)=(1,0)$, $(p,q)=(2,0)$, $(p,q)=(0,-2)$ and $(p,q)=(1,-2)$. Since $R^1(\Psi\Phi)(A)=R^2(\Psi\Phi)(A)=0$, we get $R^1\Psi R^0\Phi(A)=R^2\Psi R^0\Phi(A)=0$. If $A$ is $\D$-reflexive, $R^{-2}(\Psi\Phi)(A)=R^{-1}(\Psi\Phi)(A)=0$ and $R^0(\Psi\Phi)(A)\cong A$; hence using the edge homomorphisms, one gets 
\begin{enumerate}
\item $R^0\Psi R^{2}\Phi(A)=R^1\Psi R^{2}\Phi(A)=0$;
\item there exist the following exact sequences with natural maps
\[0\to R^0\Psi R^{1}\Phi(A)\to R^2\Psi R^{2}\Phi(A)\to A\to A/ {}_{II}E_{\infty}^{2,-2}\to 0\]
\[0\to R^1\Psi R^{1}\Phi(A)\to A/ {}_{II}E_{\infty}^{2,-2}\to R^0\Psi R^0\Phi(A)\to R^2\Psi R^{1}\Phi(A)\to 0\]
\end{enumerate}

\begin{EX}\label{ex:cinqueuno}
Let $\Lambda$ denote the $k$-algebra given by the quiver
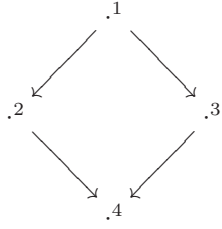
\[
\xymatrix{& {\cdot}^{1}\ar[dr]\ar[dl]\\
{\cdot}^{2} \ar[dr] && {\cdot}^{3}\ar[dl]\\
&{\cdot}^{4}}
\]
with relations such that the left projective modules are $\begin{smallmatrix} & 1 & \\ 2 & & 3\end{smallmatrix}$, $\begin{smallmatrix} 2\\4\end{smallmatrix}$, $\begin{smallmatrix} 3\\4\end{smallmatrix}$, and $\begin{smallmatrix} 4\end{smallmatrix}$. Let us consider the regular bimodule $_{\Lambda}\Lambda_{\Lambda}$; it is easy to verify that it is a cotilting bimodule of projective dimension 2. Consider the $\D$-reflexive left $\Lambda$-module $A:=\begin{smallmatrix} 2&  &3 \\  &4 & \end{smallmatrix} \oplus \begin{smallmatrix} 1 \end{smallmatrix}$. The second spectral sequence  at the second stage ${}_{II}E_2^{pq}$ and at the third and stable stage ${}_{II}E_3^{pq}={}_{II}E_{\infty}^{pq}$ assumes the following aspects:
\[
\xymatrix{
{\begin{smallmatrix} &1&\\ 2&& 3 \end{smallmatrix} \oplus \begin{smallmatrix} &1&\\ 2&& 3 \end{smallmatrix}}\ar[rrd]&0&0\\
{\begin{smallmatrix} &1&\\ 2&& 3 \end{smallmatrix} }\ar[rrd]&{\begin{smallmatrix} 4\end{smallmatrix}}&{\begin{smallmatrix} 1 \\  2 \end{smallmatrix} \oplus \begin{smallmatrix} 1\\ 3 \end{smallmatrix}}\\
0&0&{\begin{smallmatrix} 1 \\  2 \end{smallmatrix} \oplus \begin{smallmatrix} 1\\ 3 \end{smallmatrix}}
}
\qquad\qquad\qquad
\xymatrix{
{\begin{smallmatrix} 2 \end{smallmatrix} \oplus \begin{smallmatrix} 3 \end{smallmatrix}}&0&0\\
0&{\begin{smallmatrix} 4\end{smallmatrix}}&0\\
0&0&{\begin{smallmatrix} 1 \end{smallmatrix} }
}
\]

Therefore we get the following exact sequences (see the previous condition (2) in the case $n=2$):
\[0\to {\begin{smallmatrix} &1&\\ 2&& 3 \end{smallmatrix} }\to {\begin{smallmatrix} 1 \\  2 \end{smallmatrix} \oplus \begin{smallmatrix} 1\\ 3 \end{smallmatrix}}\to A= {\begin{smallmatrix} 1 \end{smallmatrix} \oplus \begin{smallmatrix} 2&&3\\ &4& \end{smallmatrix}}\to {\begin{smallmatrix} 2&&3\\ &4& \end{smallmatrix}}\to 0\]
\[0\to \begin{smallmatrix} 4 \end{smallmatrix}\to {\begin{smallmatrix} 2&&3\\ &4& \end{smallmatrix}}\to {\begin{smallmatrix} &1&\\ 2&& 3 \end{smallmatrix} \oplus \begin{smallmatrix} &1&\\ 2&& 3 \end{smallmatrix}}\to {\begin{smallmatrix} 1 \\  2 \end{smallmatrix} \oplus \begin{smallmatrix} 1\\ 3 \end{smallmatrix}}\to 0\]
\end{EX}
Note that, passing from $n=1$ to $n>1$, the spectral sequence ${}_{II}E_2^{pq}$ stabilizes at the $n+1^{th}$ stage; therefore we loose in general the possibility to describe the $\D$-reflexivity of an object $A$ in terms of properties of the objects $R^i\Psi R^j\Phi(A)$.

Resuming, the key properties which consent us to give a ``nice'' Cotilting Theorem in the 1-dimensional case are: 
\begin{description}\label{conditionI}
\item[Condition I] the spectral sequence ${}_{II}E_2^{pq}$ stabilizes at the second stage, 
\item[Condition II] the cohomologies of a $\D$-reflexive complex are $\D$-reflexive.
\end{description}
Both these properties are in general false (see Examples~\ref{ex:cinqueuno}, \ref{ex:GGdim}). 
The technical conditions assumed in the following theorems guarantee both Condition I and II. First, generalizing Theorem~\ref{cor:BB} we have

\begin{thm}\label{thm:last}
Assume $\Phi$ and $\Psi$ have cohomological dimension $n$ and 1, respectively.
An object $A$ in $\A$ is $\D$-reflexive if and only if 
\begin{enumerate}
\item
$\Phi(A)$ and $ R^1\Phi (A)$ are $\D$-reflexive;
\item $R^i\Psi R^j\Phi (A)=0$ for each $i\not= j$;
\item there exists a short exact sequence
\[0\to  R^1\Psi R^1\Phi (A)\to A\to \Psi\Phi(A)\to 0.\]
\end{enumerate}
In such a case $ R^i\Phi (A)=0$ for each $i>1$.
\end{thm}
\begin{proof}
The spectral sequence $_{II}E_2^{pq}$ stabilizes at the second stage: for only two columns survive. Therefore if $A$ is $\D$-reflexive, we get immediately the orthogonal relations
$R^i\Psi R^j\Phi (A)=0$ for each $i\not=j$.
The filtration of $A\cong R^0(\Psi\Phi(A))$ produces the short exact sequence
\[0\to R^1\Psi R^1\Phi (A)\to A\to R^0\Psi R^0\Phi(A)= \Psi\Phi(A)\to 0.\]
By the adjunction, also $\R\Phi(A)$ is $\D$-reflexive; then, by Propositions~\ref{thm:atmostone} and \ref{rem:dim}, its cohomologies $R^i\Phi(A)$ are also $\D$-reflexive. Then $\R\Psi(R^i\Phi(A))$ are $\D$-reflexive: remembering the orthogonal relations, both the complexes $\R\Psi(\Phi(A))=\Psi\Phi(A)$ and $\R\Psi(R^1\Phi(A)){\cong}R^1\Psi R^1\Phi(A)[-1]$ are $\D$-reflexive; moreover, since $0{\cong}\R\Psi(R^i\Phi(A))$,  for $i\geq 2$ the objects $R^i\Phi(A)$ are equal to zero.
Conversely, consider the triangle associated to the short exact sequence (3). By (1), (2) and the adjunction, both the complexes $\R\Psi(\Phi(A))=\Psi\Phi(A)$ and $\R\Psi(R^1\Phi(A)){\cong}R^1\Psi R^1\Phi(A)[-1]$ are $\D$-reflexive; thus one gets the $\D$-reflexivity of $A$ from the $\D$-reflexivity of the other two terms in the sequence of 3).
\end{proof}

Let us give an example where the previous theorem applies.

\begin{EX}
Let $\Lambda$ denote the $k$-algebra given by the quiver
\[
\xymatrix{{\cdot}^{0}\ar[r]& {\cdot}^{1}\ar[r]&{\cdot}^{2}\ar[r]&{\cdot}^{3}\ar[r]&{\cdot}^{4}}
\]
with relations such that the left projective modules are $\begin{smallmatrix} 0 \\ 1 \\ 2\end{smallmatrix}$, $\begin{smallmatrix} 1 \\ 2 \\ 3\end{smallmatrix}$, $\begin{smallmatrix} 2 \\ 3 \\ 4\end{smallmatrix}$, $\begin{smallmatrix} 3 \\ 4\end{smallmatrix}$ and $\begin{smallmatrix} 4\end{smallmatrix}$.
Consider the left $\Lambda$-module $_{\Lambda}U=\begin{smallmatrix} 2 \\ 3 \\ 4\end{smallmatrix}\oplus \begin{smallmatrix} 3\end{smallmatrix}\oplus \begin{smallmatrix} 1 \\ 2 \\ 3\end{smallmatrix}\oplus \begin{smallmatrix} 1\end{smallmatrix}$; it is easy to verify that it has injective dimension 2. The endomorphism ring $S:=\End ({}_{\Lambda}U)$ is the $k$-algebra given by the quiver
\[
\xymatrix{&&{\cdot}^{7}\\
{\cdot}^{5}\ar[r]&{\cdot}^{6}\ar[r]\ar[ru]&{\cdot}^{8}}
\]
with relations such that the right projective modules are $\begin{smallmatrix} 7\\ 6\end{smallmatrix}$, $\begin{smallmatrix} 8\\ 6 \end{smallmatrix}$, $\begin{smallmatrix} 6\\ 5\end{smallmatrix}$, and $\begin{smallmatrix} 5\end{smallmatrix}$. The right $S$-module $U_S=\begin{smallmatrix} 7\end{smallmatrix}\oplus \begin{smallmatrix} 7&&8\\
&6\end{smallmatrix}\oplus \begin{smallmatrix} 7\\ 6\end{smallmatrix}\oplus \begin{smallmatrix} 6\\ 5\end{smallmatrix}$ has injective dimension 1. It is easy to verify that $_{\Lambda}U_S$ is a partial cotilting bimodule. The projective module $\begin{smallmatrix} 0 \\ 1 \\ 2\end{smallmatrix}$ and its projections $\begin{smallmatrix} 0 \\ 1 \end{smallmatrix}$, $\begin{smallmatrix} 0\end{smallmatrix}$ are the only not $\D$-reflexive indecomposable $\Lambda$-modules, while all the indecomposable $S$-modules are $\D$-reflexive. In particular, consider the indecomposable left $\Lambda$-module $\begin{smallmatrix} 1\\ 2\end{smallmatrix}$; it is a $\D$-reflexive module of projective dimension 2. It satisfies the three conditions of Theorem~\ref{thm:last}:
\begin{enumerate}
\item $\Phi(\begin{smallmatrix} 1\\ 2\end{smallmatrix})=\begin{smallmatrix} 5\end{smallmatrix}$ and $R^1\Phi(\begin{smallmatrix} 1\\ 2\end{smallmatrix})=8$ are $\D$-reflexive;
\item $R^i\Psi R^j\Phi(\begin{smallmatrix} 1\\ 2\end{smallmatrix})=0$ for each $i\not=j$;
\item there exists a short exact sequence
\[0\to  R^1\Psi R^1\Phi (\begin{smallmatrix} 1\\ 2\end{smallmatrix})=\begin{smallmatrix} 2\end{smallmatrix}\to \begin{smallmatrix} 1\\ 2\end{smallmatrix}\to \begin{smallmatrix} 1\end{smallmatrix}=\Psi\Phi(\begin{smallmatrix} 1\\ 2\end{smallmatrix})\to 0.\]
\end{enumerate}
Moreover $ R^2\Phi (\begin{smallmatrix} 1\\ 2\end{smallmatrix})=0$.
\end{EX} 

A second possibility to obtain partial results is to characterize the $\D$-reflexive objects inside a suitable subclass of $\A$. Both Theorems~\ref{cor:BB} and \ref{thm:last} suggest to consider the subclass of objects $A$ in $\A$ such that $R^i\Psi R^j\Phi(A)=0$ for each $i\not=j$; in such a way, the spectral sequence  $_{II}E_2^{pq}$ stabilizes at the second stage. To have also that the cohomologies of a $\D$-reflexive complex are still $\D$-reflexive, Lemma~\ref{prop:lastt} suggests to restrict further our class (compare with \cite[Theorem~2.7]{AT}).
\begin{thm}\label{thm:lastt}
Let $(\Phi,\Psi)$ be an adjoint pair of contravariant functors of cohomological dimensions $\leq n$.
An object $A$ in $\A$, such that $ R^i\Psi R^j\Phi(A)=0$ and 
$R^i\Phi R^j\Psi R^j\Phi(A)=0$ if $i\not=j$, is $\D$-reflexive if and only if for $i=0, 1, ..., n$
\begin{enumerate}
\item $R^i\Phi(A)$ are $\D$-reflexive;
\item there exists a filtration
\[0=A_{n+1}\leq A_{n}\leq ...\leq A_0=A\]
such that $A_{i}/A_{i+1}\cong R^i\Psi R^i\Phi (A)$.
\end{enumerate}
In such a case the objects $R^i\Psi R^i\Phi (A)$ result to be $\D$-reflexive.
\end{thm}
\begin{proof}
If $A$ is $\D$-reflexive, also $\R\Phi (A)$ is $\D$-reflexive; by Lemma~\ref{prop:lastt}, {\color{black}choosing as $\rho$ the identity function}, the cohomologies $R^i\Phi(A)$, $i=1,..., n$, are $\D$-reflexive too. Since $ R^i\Psi R^j\Phi(A)=0$ for each $i\not=j$, the spectral sequence $_{II}E_2^{pq}$ stabilizes at the second stage. Therefore $R^s(\Psi\Phi)A=0$ for each $s\not=0$ and $R^0(\Psi\Phi)A\cong A$ has a finite filtration
\[0=A_{n+1}\leq A_{n}\leq ...\leq A_0=A
\]
with factors $A_{i}/A_{i+1}\cong {}_{II}E_\infty^{i, -i}={}_{II}E_2^{i, -i}= R^i\Psi R^i\Phi (A)$.

Conversely, let us assume condition (1) and (2) are satisfied. Since $R^i\Phi (A)$ is $\D$-reflexive, the complex $\R\Psi R^i\Phi(A)$ is $\D$-reflexive; we want to prove that its cohomology $R^i\Psi R^i\Phi(A)$ is $\D$-reflexive too. By hypotheses $R^j\Phi R^i\Psi R^i\Phi(A)=0$ for any $j\not=i$; moreover, since
\[R^i\Phi A\cong \R\Phi\R\Psi(R^i\Phi A)=\R\Phi(R^i\Psi R^i\Phi A[-i])=\R\Phi(R^i\Psi R^i\Phi A)[i],\]
we have
\[R^i\Phi R^i\Psi R^i\Phi A=H^{0}(\R\Phi(R^i\Psi R^i\Phi A)[i])\cong R^i\Phi A.\]
Therefore $R^j\Psi R^i\Phi R^i\Psi (R^i\Phi A)=0$ for each $j\not=i$ and by Lemma~\ref{prop:lastt} the cohomologies of $\R\Psi R^i\Phi(A)$ are $\D$-reflexive. Consider now the triangles
\[R^n\Psi R^n\Phi (A)=A_{n}\to A_{n-1}\to A_{n-1}/A_{n}=R^{n-1}\Psi R^{n-1}\Phi (A)\to A_{n}[1]\]
\[.......\]
\[A_{1}\to A_{0}=A\to A_{0}/A_{1}=R^{0}\Psi R^{0}\Phi (A)\to A_{1}[1]\]
Since $R^n\Psi R^n\Phi (A)$ and $R^{n-1}\Psi R^{n-1}\Phi (A)$ are $\D$-reflexive, also $A_{n-1}$ is $\D$-reflexive. Iterating this procedure on the other triangles, using the $\D$-reflexivity of $A_{i-1}/A_i$, $i=n,n-1,..., 1$, we prove the $\D$-reflexivity of $A$.
\end{proof}

{\color{black}We are now ready to give a Cotilting Theorem in the sense of \cite[Ch.~5]{CbF1}, between the classes of $\D$-reflexive objects induced by the pair of adjoint functors $(\Phi,\Psi)$ in the $n$-dimensional case.}

\begin{cor}\label{cor:nCotiltingTh}
{\color{black}Consider the following subclasses of the classes $\D_\A$ and $\D_\B$ of $\D$-reflexive objects in $\A$ and $\B$:
\[\overline{\D_\A}=\cap_{i\not=j}\left(\Ker R^i\Psi R^j\Phi\cap
\Ker R^i\Phi R^j\Psi R^j\Phi\right) \text{and}\]
\[\overline{\D_\B}=\cap_{i\not=j}\left(\Ker R^i\Phi R^j\Psi\cap
\Ker R^i\Psi R^j\Phi R^j\Psi\right).\]
Then setting
\[\mathcal E^i_\Phi=(\cap_{j\not=i}\Ker R^j\Phi)\cap\overline{\D_\A}\text{ and }
\mathcal E^i_\Psi=(\cap_{j\not=i}\Ker R^j\Psi)\cap\overline{\D_\B}
\]
the following conditions are satisfied:
\begin{enumerate}
\item $R^i\Phi:\overline{\D_\A}\to \mathcal E^i_\Psi$ and
$R^i\Psi:\overline{\D_\B}\to \mathcal E^i_\Phi$.
\item for each object $A$ in $\overline{\D_\A}$ and $B$ in $\overline{\D_\B}$ there exists filtrations
\[0=A_{n+1}\leq A_{n}\leq ...\leq A_0=A \text{ and}\]
\[0=B_{n+1}\leq B_{n}\leq ...\leq B_0=B\]
such that $A_{i}/A_{i+1}\cong R^i\Psi R^i\Phi (A)$ and $B_{i}/B_{i+1}\cong R^i\Phi R^i\Psi (B)$.
\item the restrictions
\[R^i\Phi:\mathcal E^i_\Phi\dual \mathcal E^i_\Psi:R^i\Psi\]
define category equivalences.
\end{enumerate}
}\end{cor}

\begin{EX}
Let $\Lambda$ denote the $k$-algebra given by the quiver ${\mathbb A}_8$
with relations such that the left projective modules are $\begin{smallmatrix}  1  \\ 2 \\ 3\\ 4\end{smallmatrix}$, $\begin{smallmatrix} 2\\3\\4\end{smallmatrix}$, $\begin{smallmatrix} 3\\4\end{smallmatrix}$, $\begin{smallmatrix} 4\\5\\6\\7\\8\end{smallmatrix}$, $\begin{smallmatrix} 5\\6\\7\\8\end{smallmatrix}$, $\begin{smallmatrix} 6\\7\\8\end{smallmatrix}$, $\begin{smallmatrix} 7\\8\end{smallmatrix}$, $\begin{smallmatrix} 8\end{smallmatrix}$. Consider the cotilting module $_{\Lambda}U= \begin{smallmatrix} 1 \\ 2\\ 3\\ 4\end{smallmatrix}\oplus\begin{smallmatrix} 1\end{smallmatrix}\oplus\begin{smallmatrix} 3\\4\end{smallmatrix}\oplus\begin{smallmatrix} 4\\5\\6\\7\\8\end{smallmatrix}\oplus\begin{smallmatrix} 5\\6\\7\\8\end{smallmatrix}\oplus\begin{smallmatrix} 6\\7\\8\end{smallmatrix}\oplus\begin{smallmatrix} 7\\8\end{smallmatrix}\oplus\begin{smallmatrix} 7\end{smallmatrix}$  of injective dimension $2$  and let $S=\End_{\Lambda}(U)$. Applying Theorem~\ref{thm:artin2} we get that any complex in $\D^b(\Lambda\lmod)$  and in $\D^b(\rmod S)$ is $\D$-reflexive w.r.t the adjunction $(\R\Hom_{\Lambda}(-, U), \R\Hom_S(-, U))$. The $\Lambda$-module $X=\begin{smallmatrix}  1  \\ 2 \\ 3\end{smallmatrix}$ satisfies the assumptions of Theorem~\ref{thm:lastt}, indeed: 
\begin{itemize}
\item $\Ext^i_S(\Ext^j_{\Lambda}(X, U), U)=0$ for $i\neq j$ 
\item  $\Hom_S(\Hom_{\Lambda}(X, U), U)=1$, $\Ext^1_S(\Ext^1_{\Lambda}(X, U), U)=2$, $\Ext^2_S(\Ext^2_{\Lambda}(X, U), U)=3$
\item $\Ext^1_{\Lambda}(1, U)=\Ext^2_{\Lambda}(1, U)=0$,  $\Hom_{\Lambda}(2, U)=\Ext^2_{\Lambda}(2, U)=0$, $\Hom_{\Lambda}(3, U)=\Ext^1_{\Lambda}(3, U)=0$
\end{itemize}
We conclude that $X$ admits a filtration $0\leq X_1 \leq X_0 \leq X$, where $X_1=3$, $X_0/X_1=2$, $X/X_0=1$.
Consider now the simple module $4$; since $\Ext^2_S(\Ext^1_{\Lambda}(4, U), U)=3$,  the assumptions of  Theorem~\ref{thm:lastt} fail.  Moreover $\Hom_S(\Hom_{\Lambda}(4, U), U)=\begin{smallmatrix} 3\\4\end{smallmatrix}$, $\Ext^1_S(\Ext^1_{\Lambda}(4, U), U)=0$, $\Ext^2_S(\Ext^2_{\Lambda}(4, U), U)=0$ and so  $4$ does not admit  any filtration with $\D$-reflexive factors.
\end{EX}

\begin{rem}
 If $R$ is noetherian and $_RU_S$ is a finitely generated cotilting bimodule, then any finitely generated projective module is reflexive and $\Hom(-, U)$-acyclic, so $\D$-reflexive. It follows that any finitely generated module is $\D$-reflexive. Let now $M\in R\lmod$ such that $\Ext^j_R(M, U)=0$ for $i\neq j$. Then, for $M$ and $\Ext^i_R(M,U)$ are $\D$-reflexive, it follows that $\Ext^j_S(\Ext^i_R(M, U))=0$  and
 \[\Ext^j_R(\Ext^i_S(\Ext^i_R(M, U),U),U)=0\text{\quad for $i\neq j$,}\] so we are in the assumption of Theorem~\ref{thm:lastt}. Thus we get {\color{black} from Corollary~\ref{cor:nCotiltingTh}} that \[M\cong \Ext^i(\Ext^i_R(M, U),U),\] recovering the Miyashita result (cfr. \cite[Theorem~1.16]{M}).
\end{rem}
\section*{Acknowledgements}
We wish to thank the referee for her/his useful suggestions which improved the paper.

\end{document}